\documentclass[12pt]{article}

\usepackage{amssymb, amsmath,amscd}
\usepackage{graphics,euscript}

\newtheorem{theorem}      {Theorem}
\newtheorem{leg}      {Legendre's Theorem}
  
\newtheorem{thm}[theorem]{Theorem}     

\newtheorem{lemma}        {Lemma}

\newtheorem{proposition}  {Proposition}

\numberwithin{corollary}{theorem}

\newcommand{\qed}{$\blacksquare$}
\newcommand{\dia}{$\diamond$}

\newtheorem{claim}       {Claim}

\newtheorem{example}      {Example}

\newtheorem{rem}       {Remark}
\newtheorem{definition}   {Definition}
\newtheorem{dfn}[definition]  {Definition}

\newcommand{\revision}[2]
    {
    \ifcase#1#2\or{\bf{#2}}\or{\bf{#2}}
    \else#2\marginpar{\small{Bad revision number!}}\fi
    }
\newcommand{\remove}[1]{}

\newcommand{\vol}{\mathrm{Vol}}
\newcommand{\cvx}{\mathrm{Conv}}
\newcommand{\diag}{\mathrm{Diag}}

\newcommand{\supp}{\mathrm{Supp}}
\newcommand{\linear}{\mathrm{Lin}}
\newcommand{\punct}[1] {\hspace{0.3cm}\text{#1}}
\newcommand{\toric}[1] {\ensuremath{\mathcal T^{#1}}}
\newcommand{\defeq}{\stackrel{\scriptstyle\mathrm{def}}{=}}
\newcommand{\romone}{\scriptscriptstyle\mathrm{I}}
\newcommand{\romtwo}{\scriptscriptstyle\mathrm{I\!I}}
\newcommand{\romthree}{\scriptscriptstyle\mathrm{I\!I\!I}}

\newcommand{\mymu}{\ensuremath{\boldsymbol{\mu}}}

\newcommand{\eps}{\varepsilon} 
\newcommand{\cF}{\mathcal{F}} 
\newcommand{\cT}{\mathcal{T}} 
\newcommand{\fA}{\EuScript{A}} 
\newcommand{\Pro}{\mathbb{P}} 
\newcommand{\pF}{{\Pro(\cF_\C(\fA))}} 
\newcommand{\Z}{\mathbb{Z}} 
\newcommand{\Zn}{\Z^n} 
\newcommand{\bO}{\mathbf{O}} 
\newcommand{\R}{\mathbb{R}} 
\newcommand{\Rn}{\R^n} 
\newcommand{\C}{\mathbb{C}} 
\newcommand{\Cs}{\C^*} 
\newcommand{\Csn}{(\Cs)^n} 
\newcommand{\Cn}{\C^n} 

\newcounter{formula}
\newlength{\formulaparindent}
\numberwithin{formula}{subsection}

\newlength{\formulawidth}
\newcommand{\formula}[2]{
\ 
\refstepcounter{formula}
\vspace{1ex} \\
\fbox{\parbox{\formulawidth}{\vspace{2ex}\hspace{\formulaparindent} Formula \theformula: #1
\[ {#2} \] }}
}

\numberwithin{equation}{subsection}

\newcount\n
\newcount\m
\def\twodigits#1{\ifnum #1<10 0\fi \number#1}

\def\hours{\n=\time \divide\n 60
  \m=-\n \multiply\m 60 \advance\m \time
    \twodigits\n :\twodigits\m}

\author{Gregorio Malajovich\thanks{Departamento de Matem\'atica Aplicada,
Universidade Federal do Rio de Janeiro, Caixa Postal 68530, CEP 21945-970,
Rio de Janeiro, RJ, Brasil. {\tt http://www.labma.ufrj.br/\~{}gregorio} ~. 
{\bf e-mail:} {\tt gregorio@labma.ufrj.br} ~.  }
\thanks{Partially supported
by CERG grants ~9040393-730, 9040402-730, and~9040188 }  
\and J. Maurice Rojas\thanks{Department of Mathematics, Texas A\&M University, 
TAMU 3368, College Station, Texas 77843-3368, USA. {\tt 
http://www.math.tamu.edu/\~{}rojas} ~, {\bf e-mail:} 
{\tt rojas@math.tamu.edu} . }
\,\thanks{Partially supported by Hong Kong UGC grant 
\#9040402-730, Hong Kong/France PROCORE Grant \#9050140-730, and a 
US National Science Foundation Mathematical Sciences Postdoctoral Fellowship.}
}

\title{
\mbox{}\\
\vspace{-3cm}
High Probability Analysis of the Condition Number of Sparse Polynomial Systems  
} 

\date{November 11, 2002}

\begin{document}

\maketitle

\thispagestyle{empty}

\begin{abstract}
Let $F\!:=\!(f_1,\ldots,f_n)$ be a random polynomial system with 
fixed $n$-tuple of supports. Our main result is an upper bound on the 
probability that 
the condition number of $f$ in a region $U$ is larger than $1/\eps$. The bound 
depends on an integral of a differential form on a toric manifold and 
admits a simple explicit upper bound when the Newton polytopes (and underlying 
covariances) are all identical.  

We also consider polynomials with real coefficients and give bounds
for the expected number of real roots and (restricted) condition number.  
Using a K\"ahler geometric framework throughout, we also express the expected 
number of roots of $f$ inside a region $U$ as the integral over $U$ of a 
certain {\bf mixed volume}  form, thus recovering the classical mixed volume 
when $U = (\mathbb C^*)^n$. 
\end{abstract}


\noindent 
{\bf 2000 Math Subject Classification:} 
65H10, 
52A39. 

\section{Introduction}
From the point of view of numerical analysis, it is not only the number 
of complex solutions of a polynomial system which make it hard to solve 
numerically but the sensitivity of its roots to small perturbations in the 
coefficients. This is formalized in the condition number, $\mu(f,\zeta)$ 
(cf.\ Definition \ref{condnum} of Section \ref{sub:stronger}), which dates back to work of 
Alan Turing \cite{TURING}. 
In essence, $\mu(f,\zeta)$ measures the sensitivity of a solution $\zeta$ to perturbations in 
a problem $f$, and a large condition number is meant to imply that $f$ is intrinsically hard 
to solve numerically. Such analysis of numerical conditioning, while having been applied for 
decades in numerical linear algebra (see, e.g., \cite{DEMMEL}), has only been applied to 
computational algebraic geometry toward the end of the twentieth century (see, e.g., 
\cite{BEZ1}). 

Here we use K\"ahler geometry to analyze the numerical conditioning of sparse polynomial 
systems, thus setting the stage for more realistic complexity bounds for the numerical 
solution of polynomial systems. Our bounds generalize some earlier results of 
Kostlan~\cite{KOSTLAN} and Shub and Smale~\cite{BEZ4} on the more restricted dense case, 
and also yield new formulae for the expected number of roots (real and complex) in a region. 
The appellations ``sparse'' and ``dense'' respectively refer to either (a) taking into account 
the underlying monomial term structure or (b) ignoring this finer structure and simply 
working with degrees of polynomials. Since many polynomial systems occuring in practice 
have rather restricted monomial term structure, sparsity is an important consideration 
and we therefore strive to state our complexity bounds in terms of this refined information.  

To give the flavor of our results, let us first make some necessary 
definitions. We must first formalize the spaces of polynomial systems 
we work with and how we measure perturbations in the spaces of 
problems and solutions. 

\begin{dfn}
Given any finite subset $A\!\subset\Zn$, let $\cF_\C(A)$ (resp.\ 
$\cF_\R(A)$) denote the vector space of all polynomials in $\C[x_1,\ldots,x_n]$  
(resp.\ $\R[x_1,\ldots,x_n]$) of the form 
$\sum\limits_{a \in A} c_a x^a$ where the notation $x^a\!:=\!x^{a_1}\cdots x^{a_n}$ 
is understood. For any finite subsets $A_1,\ldots,A_n\!\subset\!\Zn$ we then 
let $\fA:=(A_1,\ldots,A_n)$ and $\cF_\C(\fA)\!:=\!\cF_\C(A_1)\times \cdots 
\times \cF_\C(A_n)$ (resp.\ $\cF_\R(\fA)\!:=\!\cF_\R(A_1)\times \cdots \times 
\cF_\R(A_n)$). \dia 
\end{dfn} 
The $n$-tuple $\fA$ will thus govern our notion of sparsity as well as the perturbations 
allowed in the coefficients of our polynomial systems. It is then easy 
to speak of random polynomial systems and the distance to the nearest degenerate system. 
Recall that a {\bf degenerate} root of $f$ is simply a root of $f$ having Jacobian of 
rank $<\!n$. 
\begin{dfn} 
By a {\bf complex (resp.\ real) random sparse polynomial system} we will mean 
a choice of $\fA\!:=\!(A_1, \ldots,A_n)$ and an assignment of a probability 
measure to each $\cF_\C(A_i)$ (resp.\ $\cF_\R(A_i)$) as follows: endow 
$\cF_\C(A_i)$ (resp.\ $\cF_\R(A_i)$) with an independent complex (resp.\ real) 
Gaussian distribution having mean $\bO$ and a (positive definite and diagonal) 
covariance matrix $C_i$. Finally, let the {\bf discriminant variety}, 
$\Sigma(\fA)$, denote the set of all $f\!\in\!\cF_\C(\fA)$ (resp.\ 
$f\!\in\!\cF_\R(\fA)$) with a degenerate root and define 
$\cF_\zeta(\fA)\!:=\!\{f\!\in\!  \cF_\C(\fA) \; | \; f(\zeta)\!=\!\bO\}$ 
(resp.\ $\cF_\zeta(\fA)\!:=\!\{f\!\in\!  \cF_\R(\fA) \; | \; 
f(\zeta)\!=\!\bO\}$) and $\Sigma_{\zeta}(\fA)\!:=\!\cF_\zeta(\fA)\cap 
\Sigma(\fA)$. \dia 
\end{dfn} 

\begin{thm}
\label{thm:simpleunmixed}
Suppose $A\!\subset\!\Zn$ is a finite set with a convex hull of 
positive volume and $\fA\!:=\!(\underset{n}{\underbrace{A,\ldots,A}})$. 
Then there is a natural metric $d(\cdot,\cdot)$ on $\cF_\C(\fA)$ 
such that $\mu(f,\zeta)\!=\!\frac{1}{d(f,\Sigma_{\zeta}(\fA))}$. Furthermore,  
\[ \mathrm{Prob}\left[\mu(f,\zeta)\!\geq\!
\frac{1}{\eps} \text{ for some root } \zeta\!\in\!\Csn \text{ of } f 
\right]\leq 
n^3(n+1)\vol(A)(\#A-1)(\#A-2)\eps^4 \]
where $f$ is a complex random sparse polynomial system, $\#A$ denotes the 
number of points in $A$, and $\vol(A)$ denotes 
the volume of the convex hull of $A$ (normalized so that 
$\vol(\bO,e_1,\ldots,e_n)\!=\!1$). 
\end{thm} 

\noindent
The above theorem is in fact a simple corollary of two much more general 
theorems (Theorems \ref{thm:condnumber} and \ref{mixed-cond}) which also include as a special 
case an analogous result of Shub and Smale in the dense case \cite[Thm.\ 1, Pg.\ 237]{BCSS}. 
We also note that theorems such as the one above are natural precursors to explicit bounds 
on the number of steps required for a homotopy algorithm \cite{BEZ1} to solve $f$. 
We will pursue the latter topic in a future paper. Indeed, one of our long term goals is to 
provide a rigourous and explicit complexity analysis of the numerical homotopy algorithms for 
sparse polynomial systems developed by Verschelde et.\ al.~\cite{VVC}, 
Huber and Sturmfels~\cite{HS}, and Li and Li~\cite{LI}. 

The framework underlying our first main theorem involves K\"ahler geometry, 
which is the intersection of Riemannian metrics and symplectic and complex 
structures on manifolds. On a more concrete level, we can give new formulae for the 
expected number of roots of $f$ in a region $U$. For technical reasons, we will mainly 
work with {\bf logarithmic coordinates}. That is, we will let 
$\cT^n$ be the $n$-fold product of cylinders $(\R\times (\R \ \mathrm{mod} \ 
2\pi))^n\!\subset\!\Cn$,  
and use coordinates $p+iq\!:=\!(p_1+iq_1,\ldots,p_n+iq_n)\!\in\!\cT^n$ to stand for a root 
$\zeta\!:=\!\exp(p+iq)\!:=\!(e^{p_1+iq_1},\ldots,e^{p_n+iq_n})$ of $f$.  
Roots with zero coordinates can be handled by then working in a suitable 
{\bf toric compactification} and this is made precise in Section \ref{symplectic}. 
The idea of working with roots of polynomial systems
in logarithmic coordinates seems to be extremely 
classical, yet it gives rise to interesting and
surprising connections (see the discussions in~\cite{GRAEFFE, 
TANGRA,VIRO}).
\begin{theorem}\label{GEN-BERNSHTEIN}
Let $A_1,\ldots,A_n$ be finite subsets of $\Zn$ and  
$U\!\subseteq\!\cT^n$ be a measurable region. Pick positive definite 
diagonal covariance matrices $C_1,\ldots,C_n$ 
and consider a complex random 
polynomial system specified by the data $(A_1,C_1,\ldots,A_n,C_n)$.  
Then there are natural real $2$-forms $\omega_{A_1},\ldots,
\omega_{A_n}$ on $\cT^n$ such that the expected number of roots of 
$f$ in $\exp U \subseteq (\mathbb C^*)^n$ is exactly 
\[ \frac{(-1)^{n(n-1)/2}}{\pi^n}\ \int_U 
\omega_{A_1} \wedge \cdots \wedge \omega_{A_n}
\punct{.} \] 
In particular, when $U\!=\!\Csn$, the above expression is exactly the 
mixed volume of the convex hulls of $A_1,\ldots,A_n$ (normalized so 
that the mixed volume of $n$ standard $n$-simplices is $1$).  
\end{theorem}
\noindent 
See~\cite{BUZA,SANGWINEYAGER} for the classical definition of
mixed volume and its main properties. The result above generalizes the famous 
connection between root counting and mixed volumes discovered 
by David N.\ Bernshtein~\cite{BERNSHTEIN}. The special case of unmixed systems with 
identical coefficient distributions ($A_1 = \cdots = A_n$, $C_1 = \cdots = C_n$) recovers 
a particular case of Theorem~8.1 in~\cite{EDELMAN-KOSTLAN}. However, 
comparing Theorem \ref{GEN-BERNSHTEIN} and \cite[Theorem 8.1]{EDELMAN-KOSTLAN}, 
this is the only overlap since neither theorem generalizes the other.
The very last assertion of Theorem \ref{GEN-BERNSHTEIN} (for uniform variance $C_i\!=\!I$ 
for all $i$) was certainly known to Gromov \cite{GROMOV}, and a version of 
Theorem \ref{GEN-BERNSHTEIN} was known to Kazarnovskii \cite[p.~351]{KAZARNOVSKII} and 
Khovanskii \cite[Prop.\ 1, Sec.\ 1.13]{KHOVANSKII}.  
In \cite{KAZARNOVSKII}, the supports $A_i$ are even allowed
to have complex exponents.  However, uniform variance is again assumed. 
His method may imply this special case of Theorem~\ref{GEN-BERNSHTEIN}, but the indications
given in~\cite{KAZARNOVSKII} were insufficient for us to reconstruct a proof.

As a consequence of our last result, we can also give a coarse estimate on 
the expected number of real roots in a region. 
\begin{theorem} \label{expected-real}
Let $U$ be a measurable subset of $\Rn$ 
with Lebesgue volume $\lambda(U)$. Then, following the notation above, suppose 
instead that $f$ is a real random polynomial system. Then the average number of real roots 
of $f$ in $\exp U \subset \R^n_+$ is bounded above by
\[ (4\pi^2)^{-n/2} \sqrt{\lambda(U)} \sqrt{ \int_{(p,q)\in U\times [0,2\pi)^n} 
(-1)^{n(n-1)/2} \omega_{A_1} \wedge \cdots \wedge \omega_{A_n}}. 
\]
\end{theorem}

\noindent
This bound is of interest when $n$ and $U$ are fixed, 
in which case the expected number of positive real roots grows as
the square root of the mixed volume.

\subsection{Stronger Results Via Mixed Metrics}
\label{sub:stronger} 
Our remaining new results, which further sharpen the preceding bounds 
and formulae, will require some additional notation. 
\begin{dfn} 
We define a norm on $\cF_\C(A_i)$ by $\left\| f^i \right\|_{C^{-1}}\!:=\!f^i 
C^{-1}(f^i)^H$ where, in the last expression, we consider $f^i$ via its 
row vector of coefficients and $(\cdot)^H$ denotes the usual Hermitian conjugate transpose. 
Finally, we define a norm on $\cF_\C(\fA)$ by 
$\left\|f\right\|\!:=\!\sum\limits^n_{i=1} \left\| f^i \right\|_{C^{-1}_i}$, and a 
metric $d_{\Pro}$ on the product of projective spaces $\pF\!:=\!\Pro(\cF_\C(A_1))
\times \cdots \times \Pro(\cF_\C(A_n))$ by 
$d_{\Pro}(f,g)\!:=\!\sum\limits^n_{i=1} 
\min\limits_{\lambda\in\Cs} \frac{\left\|f^i - \lambda g^i\right\|}{\left\|f^i\right\|}$,  
where we implicitly use the natural embedding of $\Pro(\cF_\C(A_i))$ into 
the unit hemisphere of $\cF_\C(A_i)$. \dia 
\end{dfn} 
Each of the terms in the sum above corresponds to the square of the
sine of the Fubini (or angular) distance between $f^i$ and $g^i$.
Therefore, $d_{\mathbb P}$ is never larger than the Hermitian distance 
between points in $\cF_\C(\fA)$, but is a correct first-order approximation of the 
distance when $g \rightarrow f$ in $\mathbb P (\cF_\C(\fA))$ (compare 
with \cite[Ch.~12]{BCSS}).

Recall that $T_pM$ denotes the tangent space at $p$ of a manifold $M$. 
\begin{dfn} 
\label{condnum}
Define the {\bf evaluation map}, $ev_\fA$, as follows: 
\[
\begin{array}{lrcl}
\mathit{ev_\fA}: & \mathcal F \times \toric{n} & \rightarrow & \mathbb C^n \\
& ((f^1,\ldots,f^n), p+iq) & \mapsto & (f^1(\exp(p+iq)),\ldots,f^n(\exp(p+iq))). 
\end{array}
\]
Given any root $\exp(p+iq)$ of an $f$ in $\cF_\C(\fA)$, the {\bf condition number of $f$ at 
$p+iq$}, $\mymu(f,p+iq)$, is then defined to be the 
operator norm 
\[ \left\|\left. DG\right|_f\right\|\!:=\!\max\limits_{\|g\|=1} \left\| 
\left. DG\right|_f\right\|,\] where $G$ is the unique branch of the implicit function which 
satisfies $G(f)\!=\!p+iq$ and $ev_\fA(g,G(g))\!=\!\bO$ for all $g$ sufficiently near $f$, and 
$DG : T_f\cF_\C(\fA) \longrightarrow T_{p+iq} \cT^n$ is the derivative of $G$. (We 
set the condition number $\mu(f,p+iq)\!:=\!+\infty$ in the event that 
$Df$ is zero and $G$ thus fails to be uniquely defined.) \dia 
\end{dfn} 
\noindent
Note that the implied norm on $T_f \cF_\C(\fA)$ was detailed in the previous definition, while 
the implied norm on $T_{p+iq} \cT^n$ has intentionally been left unspecified. 
This is because while $\cF_\C(\fA)$ admits a natural Hermitian structure, the solution-space
\toric{n} admits $n$ different natural Hermitian structures (one from each 
support $A_i$, as we shall see in the next section). Nevertheless, we can 
give useful bounds on the condition number and give an unamibiguous definition 
in certain cases. 
\begin{theorem}[Condition Number Theorem] \label{thm:condnumber}
If $(p,q)\!\in\cT^n$ is a non-degenerate root of $f$ then 
\[
\max_{ \|\dot f\| \le 1}
\min_{i}
\| DG_f \dot f \|_{A_i}
\le
\frac{1}{d_{\mathbb P}(f, \Sigma_{(p,q)})}
\le
\max_{ \|\dot f\| \le 1}
\max_{i}
\| DG_f \dot f \|_{A_i}
\punct{.}
\]
In particular, if $A_1 = \cdots = A_n$ and $C_1 = \cdots = C_n$, then
\[
\max_{ \|\dot f\| \le 1}
\min_{i}
\| DG_f \dot f \|_{A_i}
=
\max_{i}
\max_{ \|\dot f\| \le 1}
\| DG_f \dot f \|_{A_i}
=
\frac{1}{d_{\mathbb P}(f, \Sigma_{(p,q)})}
\]
and we can define $\mymu(f;(p,q))$ to be any of the three preceding 
quantities. 
\end{theorem}
This generalizes \cite[Thm.\ 3, pg.\ 234]{BCSS} which is essentially equivalent 
to the last assertion above, in the special case where $A_i$ is an $n$-column matrix 
whose rows $\{A^\alpha_i\}_\alpha$ consist of all partitions of $d_i$ into $n$ non-negative 
integers and $C_i\! = \!\diag_{\alpha} \left( \frac{(d_i-1)!}{ 
(A_i)^{\alpha}_1 !
(A_i)^{\alpha}_2 !
\cdots
(A_i)^{\alpha}_n !
(d_i - \sum_{j=1}^n (A_i)^{\alpha}_j) !
} \right)$ --- in short, the case where one considers complex random polynomial systems with 
$f^i$ a degree $d_i$ polynomial and the underlying probability measure is 
invariant under a natural action of the unitary group $U(n+1)$ on the space of roots. 
The last assertion of Theorem \ref{thm:condnumber} also 
bears some similarity to Theorem D of \cite{DEDIEU-PARK-CITY} where the notion of metric is 
considerably loosened to give a statement which applies to an even more general class 
of equations. However, our philosophy is radically different: we consider the inner
product in $\cF_\C(\fA)$ as the starting point of our investigation and we do {\bf not} 
change the metric in the fiber $\mathcal F_{(p,q)}$. Theorem \ref{thm:condnumber} 
thus gives us some insight about reasonable intrinsic metric structures on $\mathcal T^n$.

In view of the preceding theorem, we can define a restricted condition number with respect to
any measurable sub-region $U\subset\cT^n$ as follows: 
\begin{dfn} 
\label{intrinsic}
We let $\mymu (f; U)\!:=\!\frac{1}{\min_{(p,q) \in U} d_{\mathbb P}(f, \Sigma_{(p,q)})}$. 
Also, via the natural $GL(n)$-action on $T_{(p,q)}\mathcal T^n$ defined by 
$(\dot p,\dot q) \mapsto (L\dot p,L\dot q)$ for any $L \in GL(n)$, 
we define the {\em mixed dilation} of the tuple $(\omega_{A_1}, \cdots, \omega_{A_n})$
as:
\[
\kappa(\omega_{A_1},\cdots,\omega_{A_n}; (p,q)) :=
\min_{L \in GL(n)}
\max_i
\frac{
  \max_{\|u\|=1} (\omega_{A_i})_{(p,q)} (Lu,JLu)
}{
  \min_{\|u\|=1} (\omega_{A_i})_{(p,q)} (Lu,JLu)
}
\punct{,}
\]
where $J : T \cT^n \longrightarrow T \cT^n$ is canonical complex structure of $\cT^n$. 
Finally, we define 
$\kappa_{U} := \sup_{(p,q) \in U} \kappa(\omega_{A_1},\cdots,\omega_{A_n}; (p,q))$, 
provided the supremum exists, and $\kappa_U\!:=\!+\infty$ otherwise. \dia 
\end{dfn}

We can then bound the expected number of roots with condition number
$\mymu > \eps^{-1}$ on $U$ in terms of the mixed volume form,
the mixed dilation $\kappa_U$ and the expected 
number of ill-conditioned roots in the {\em linear case}. The linear case corresponds to 
the point sets and covariance matrices below:

\begin{align*}
A^{\linear}_i &=
\left[
\begin{matrix}
0 & \cdots & 0 \\
1 & & \\
  & \ddots & \\
  & & 1  
\end{matrix}
\right]
&
C^{\linear}_i =
\left[
\begin{matrix}
1 \\
  & 1 \\
  &   & \ddots & \\
  &   &        & 1  
\end{matrix}
\right]
\end{align*} 

\begin{theorem}[Expected Value of the Condition Number]
\label{mixed-cond} 
Let $\nu^{\linear}(n,\eps)$ be the probability that a complex random 
system of $n$ polynomial in $n$ variables has condition number
larger than $\eps^{-1}$. Let $\nu^{A}(U,\eps)$ be the
probability that $\mymu(f,U) > \eps^{-1}$ for a complex random polynomial system
$f$ with supports $A_1,\cdots,A_n$ and covariance matrices $C_1,\cdots,C_n$.

Then, $\nu^A(U,\eps) \le \frac{ \int_U \bigwedge \omega_{A_i} }
{ \int_U \bigwedge \omega_{A^{\linear}_i}} \
\nu^{\linear} (n, \sqrt{\kappa_U} \eps)$. 
\end{theorem}

Our final main result concerns the distribution of the real roots 
of a real random polynomial system. Let $\nu_{\mathbb R}(n,\eps)$ be the probability that 
a real random linear system of $n$ polynomials in $n$ variables has condition number larger 
than $\eps^{-1}$.  
\begin{theorem}\label{unmixed:real}
  Let $A = A_1 = \cdots = A_n$ and $C=C_1 = \cdots = C_n$ and 
let $U \subseteq \mathbb R^n$ be measurable. Let $f$
be a real random polynomial system. Then,
\[
\mathrm{Prob} \left[ \mymu(f,U) > \eps^{-1} \right]
\le
E(U) \ \nu_{\mathbb R}(n,\eps)
\]
where $E(U)$ is the expected number of real roots on $U$.
\end{theorem}

\noindent 
Note that $E(U)$ depends on $C$, so even if we make
$U = \mathbb R^n$ we may still obtain a bound depending on $C$.
Shub and Smale showed in~\cite{BEZ2} that the expected number of real roots in the dense
case (with a particular choice of probability measure)
is exactly the square root of the expected number of complex roots.
The sparse analogue of this result seems 
hard prove even in the general unmixed case: Explicit formul\ae \ for the unmixed case
are known only in certain special cases, e.g., certain systems of bounded multi-degree 
\cite{ROJAS-PARK-CITY,MCLENNAN}. Hence our last theorem can be interpreted as another 
step toward a fuller generalization. 

\section{Symplectic Geometry and Polynomial Systems}
\label{symplectic}

\subsection{Some Basic Definitions and Examples}\label{sec:symp}
For the standard definitions and properties of symplectic structures, 
complex structures, Riemannian manifolds, and K\"ahler manifolds, we refer the 
reader to \cite{MCDUFF-SALAMON,CHERN-CHEN-LAM}. A treatment focusing on toric 
manifolds can be found in \cite[Appendix A]{GUILLEMIN}. We briefly review 
a few of the basics before moving on to the proofs of our theorems.  
\begin{definition}[K\"ahler manifolds]\label{def4}
Let $M$ be a complex manifold, with complex structure $J$ and a strictly 
positive symplectic $(1,1)$-form $\omega$ on $M$ (considered as a real manifold). 
We then call the triple $(M, \omega, J)$ a {\bf K\"ahler manifold}. \dia 
\end{definition}

\begin{example}[Affine Space] 
  We identify $\mathbb C^{M}$ with $\mathbb R^{2M}$ and use coordinates
$Z^i = X^i + \sqrt{-1} Y^i$. The {\em canonical} $2$-form 
$\omega_Z = \sum_{i=1}^M dX_i \wedge dY_i$ makes $\mathbb C^{M}$ into a
symplectic manifold.

  The natural complex structure $J$ is just the multiplication by ${\sqrt{-1}}$.
The triple $(\mathbb C^M, \omega_Z, J)$ is a K\"ahler manifold. \dia  
\end{example}

\begin{example}[Projective Space]\label{ex:proj}
  Projective space $\mathbb P^{M-1}$ admits a {\em canonical} $2$-form
defined as follows. Let $Z = (Z^1, \cdots, Z^M) \in (\mathbb C^{M})^*$,
and let $[Z] = (Z^1 : \cdots : Z^M) \in \mathbb P^{M-1}$ be the corresponding 
point in $\mathbb P^{M-1}$.
  The tangent space $T_{[Z]} \mathbb P^{M-1}$ may be modelled by $Z^{\perp}
\subset T_Z \mathbb C^M$. 
  Then we can define a two-form on $\mathbb P^{M-1}$ by setting:
\[
\omega_{[Z]} (u,v) = \|Z\|^{-2} \omega_Z (u,v) \punct{,}
\]
  where it is assumed that $u$ and $v$ are orthogonal to $Z$. The 
latter assumption tends to be quite inconvenient, and most people 
prefer to pull $\omega_{[Z]}$ back to $\mathbb C^M$ by the canonical 
projection $\pi: Z \mapsto [Z]$.
  It is standard to write the pull-back $\tau = \pi^* \omega_{[Z]}$
as:
\[
\tau_Z = - \frac{1}{2}d J^* d \ \frac{1}{2}\log \|Z\|^2 \punct{,}
\]
using the notation $d\eta = \sum_i \frac{\partial \eta}{p_i} \wedge dp_i
+ \frac{\partial \eta}{q_i} \wedge dq_i$, and where $J^*$ denotes the
pull-back by $J$.

  Projective space also inherits the complex structure from $\mathbb C^M$.
Then $\omega_{[Z]}$ is a strictly positive $(1,1)$-form. The corresponding
metric is called {\em Fubini-Study} metric in $\mathbb C^M$ or 
$\mathbb C^{M-1}$. \dia 
\end{example}

\begin{rem}
  Some authors prefer to write 
$\sqrt{-1} \partial \bar \partial$
instead of $- \frac{1}{2} d J^* d$. The following notation
is assumed:
$\partial \eta = \sum_i \frac{\partial \eta}{Z_i} \wedge dZ_i$
and $\bar \partial \eta = \sum_i \frac{\partial \eta}{\bar Z_i} \wedge d\bar Z_i$.
Then they write $\tau_Z$ as:
\[
\tau_Z = \frac{\sqrt{-1}}{2}
\left(
\frac{ \sum_i dZ_i \wedge d\bar Z_i }{\|Z\|^2}
-
\frac{ \sum_i Z_i d\bar Z_i \wedge \sum_j \bar Z_j dZ_j}{\|Z\|^4}
\right)
\punct{. \dia}
\]
\end{rem}

\begin{example}[Toric K\"ahler Manifolds from Point Sets] 
Let $A$ be any $M \times n$ matrix with integer entries whose row vectors 
have $n$-dimensional convex hull and let $C$ be any diagonal positive 
definite $n\ times n$ matrix. Define the map $\hat V_A$ from $\mathbb C^n$ into 
$\mathbb C^{M}$ by 
$\hat V_A: z \mapsto C^{1/2}
\left[
\begin{matrix} z^{A^1} \\ \vdots \\ z^{A^M} \end{matrix} \right]$. 
We can also compose with the projection into projective space to obtain 
a slightly different map 
$V_A = \pi \circ \hat V_A: \mathbb C^n \rightarrow \mathbb P^{M-1}$ defined 
by $V_A: z \mapsto [\hat V_A(z)]$.
When $C$ is the identity, the Zariski closure of the
image of $V_A$ is called the {\em Veronese variety} and the map $V_A$ is
called the {\em Veronese embedding}. Note that $V_A$ is not defined
for certain values of $z$, like $z=0$. Those values comprise the
{\em exceptional set} which is a subset of the coordinate hyperplanes. 

There is then a natural symplectic structure on the closure of the image of
$V_A$, given by the restriction of the Fubini-Study $2$-form $\tau$: 
We will see below (Lemma~\ref{dimnondeg}) 
that by our assumption on the convex hull of the rows of $A$, we have that $DV_A$ is of rank 
$n$ for $z \in (\mathbb C^*)^n$. 
Thus, by the above lemma, we can pull-back this structure to $\Csn$ by 
$\Omega_A = V_A ^* \tau$. Also, we can pull back the complex structure of $\mathbb P^{M-1}$,
so that $\Omega_A$ becomes a strictly positive $(1,1)$-form. 
Therefore, the matrix $A$ defines a K\"ahler manifold
$(\mathbb \Csn, \Omega_A, J)$. \dia   
\end{example}
\begin{lemma} \label{dimnondeg}
Let $A$ be a matrix with non-negative integer entries,
such that $\cvx(A)$ has dimension $n$. Then $(Dv_A)_p$ is injective,
for all $p \in \mathbb R^n$.
\end{lemma}

\noindent
{\bf Proof:}  The conclusion of this Lemma can fail only if there are $p \in \mathbb R^n$
and $u\ne 0$
with $(Dv_A)_p u = 0$. This means that
\[
P_{v_A(p)} \mathrm{diag} (v_A)_p A u = 0
\punct{.}
\]

This can only happen if $\mathrm{diag} (v_A)_p A u$ is in the
space spanned by $(v_A)_p$, or, equivalently, $Au$ is in the
space spanned by $(1,1, \cdots, 1)^T$. This means that all the
rows $a$ of $A$ satisfy $a u = \lambda$ for some $\lambda$.
Interpreting a row of $A$ as
a vertex of $\cvx (A)$, this means that $\cvx (A)$ is contained in
the affine plane $\{ a : a u = \lambda \}$. \qed

\begin{example}[Toric Manifolds in Logarithmic Coordinates] 
For any matrix $A$ as in the previous example, we can pull-back the K\"ahler
structure of $((\mathbb C^*)^n, \Omega_A, J)$ to obtain another K\"ahler
manifold $(\toric{n}, \omega_A, J)$.
(Actually, it is the same object in
logarithmic coordinates, minus points at ``infinity''.) 
An equivalent definition is to pull back the K\"ahler structure
of the Veronese variety by $\hat v_A \defeq \hat V_A \circ \exp$. \dia
\end{example}

\begin{rem}
  The Fubini-Study metric on $\mathbb C^M$ was constructed by applying the
operator $-\frac{1}{2}dJ^*d$ to a certain convex function
(in our case, $\frac{1}{2}\log \|Z\|^2$). This is a general standard way to construct
K\"ahler structures. In~\cite{GROMOV}, it is explained how to associate a
(non-unique) convex function to any convex body, thus producing an associated
K\"ahler metric. \dia 
\end{rem}

For the record, we state explicit formul\ae\ for several of the
invariants associated to the K\"ahler manifold
$(\toric{n}, \omega_A, J)$. First of all, the function
$g_A = g \,\circ \,\hat v_A$ is precisely:

\formula{\label{211}The canonical Integral $g_A$ (or {\em K\"ahler potential})
of the convex set associated to $A$}
{ g_A(p):=\frac{1}{2} \log 
\left(
\left( \exp (A \cdot p) \right) ^T
C
\left( \exp (A \cdot p) \right) \right)} 

The terminology {\em integral} is borrowed from mechanics,
and it refers to the invariance of $g_A$ under a $[0,2\Pi)^n$-action. Also, the
gradient of $g_A$ is called the {\em momentum map}. Recall that 
the Veronese
embedding takes values in projective space. We will use the
following notation:
$v_A(p) = \hat v_A(p) / \| \hat v_A(p) \|$. 
This is independent
of the representative of equivalence class $v_A(p)$. 
Now, let $v_A(p)^2$ mean
coordinatewise squaring and $v_A(p)^{2T}$ be the transpose 
of $v_A(p)^2$. The gradient of $g_A$ is then:
\formula{The Momentum Map associated to $A$}
{
\nabla g_A = v_A(p)^{2T} A
}

\formula{\label{secder} Second derivative of $g_A$}
{
D^2 g_A = 2 Dv_A(p)^T Dv_A(p)
}

We also have the following formulae:  
\formula{\label{omega} The symplectic $2$-form associated to $A$:}
{
(\omega_A)_{(p,q)} = \frac{1}{2}\sum_{ij} (D^2 g_A)_{ij} dp_i \wedge dq_j
}

\formula{Hermitian structure of \toric{n} associated to $A$:}
{
(\langle u, w \rangle_{\scriptscriptstyle A})_{(p,q)} 
= 
u^H (\frac{1}{2}D^2 g_A)_p w
}

In general, the function $v_A$ goes from \toric{n}
into projective space. Therefore, its derivative
is a mapping

\[
(Dv_A)_{(p,q)}:
T_{(p,q)}\toric{n}
\rightarrow T_{\scriptscriptstyle v_A(p+q \sqrt{-1})} \mathbb P^{M-1}
\simeq \hat v_A(p+q\sqrt{-1}) ^{\perp} \subset \mathbb C^M
\punct{.}
\]

For convenience, we will write this derivative as a mapping into
$\mathbb C^M$, with range $\hat v_A(p+q\sqrt{-1}) ^{\perp}$.
Let $P_v$ be the projection operator
\[
P_v = I - \frac{1}{\|v\|^2} v v^H
\punct{.}
\]
We then have the following formula. 
\formula{\label{deriv} Derivative of $v_A$}
{
(Dv_A)_{\scriptscriptstyle (p,q)} = 
P_{ \hat v_A(p+q \sqrt{-1}) }
\diag \left( \frac{ \hat v_A(p+q \sqrt{-1}) }{\|\hat v_A(p+q \sqrt{-1}\|}
\right) A
}

An immediate consequence of Formula~\ref{deriv} is:

\begin{lemma}\label{lem-orthog}
  Let $f \in \mathcal F_A$ and $(p,q) \in \toric{n}$
be such that $f \cdot \hat v_A(p+q\sqrt{-1})=0$. Then, 
$f \cdot (Dv_A)_{(p,q)} = \frac{1}{\| \hat v_A(p,q)\|} f \cdot (D\hat v_A)_{(p,q)}$
\end{lemma}

  In other words, when $(f \circ \exp) (p+q\sqrt{-1})$ vanishes,
$Dv_A$ and $D\hat v_A$ are the same up to scaling.
Noting that the Hermitian metric can be written 
$(\langle u,w\rangle_{\scriptscriptstyle A})_{(p,q)} = u^h Dv_A(p,q)^H Dv_A(p,q) w$, we 
also obtain the following formula.  
\formula{Volume element of $(\toric{n}, \omega_A, J)$}
{
d\toric{n}_A = \det \left( \frac{1}{2} \ D^2 g_A(p)\right) \  
dp_1 \wedge \cdots \wedge dp_n \wedge dq_1 \wedge \cdots \wedge dq_n 
}

\subsection{Toric Actions and the Momentum Map}

The {\em momentum map}, also called {\em moment map}, was 
introduced in its modern formulation by Smale~\cite{SMALE70}
and Souriau~\cite{SOURIAU}. The reader may consult one of
the many textbooks in the subject (such as Abraham and
Marsden~\cite{ABRAHAM-MARSDEN} or McDuff and Salamon 
\cite{MCDUFF-SALAMON}) for a general exposition. 

In this section we instead follow the point of view
of Gromov~\cite{GROMOV}. The main results in this
section are the two propositions below.  
\begin{proposition} \label{convexity}
  The momentum map $\nabla g_A$ maps 
$\toric{n}$ onto the interior of $\cvx(A)$.
When $\nabla g_A$ is restricted to the
real $n$-plane $[q=0] \subset \toric{n}$,
this mapping is a bijection.  
\end{proposition}

This would appear to be a particular case of the 
Atiyah-Guillemin-Sternberg theorem \cite{ATIYAH82,GS}. However, technical
difficulties prevent us from directly applying this result here.\footnote{The 
Atiyah-Guillemin-Sternberg applies to compact symplectic manifolds and the 
implied compactification of $\cT^n$ may have singularities.}   

\begin{proposition} \label{volpres}
  The momentum map $\nabla g_A$ is a volume-preserving
map from the manifold $(\toric{n}, \omega_A, J)$
into $\cvx(A)$, up to a constant, in the following
sense: if $U$ is a measurable region of $\cvx(A)$,
then
\[
\vol \left( (\nabla g_A)^{-1}(U) \right) = \pi^{n}\  \vol (U)
\punct{. }
\]
\end{proposition}

\noindent 
{\bf Proof of Proposition~\ref{volpres}:} 
Consider the mapping
\[
\begin{array}{rrcl}
M: & \toric{n} & \rightarrow &  \frac{1}{2}\cvx(A) \times \mathbb T^n \\
& (p,q) & \mapsto & (\frac{1}{2}\nabla g_A(p), q)
\end{array}
\punct{.}
\]

Since we assume $\dim \cvx(A) = n$, we can apply Proposition \ref{convexity}
and conclude that $M$ is a diffeomorphism. 
\medskip

The pull-back of the canonical symplectic structure in
$\mathbb R^{2n}$ by $M$ is precisely $\omega_A$, because
of Formul\ae~\ref{secder} and~\ref{omega}. 
Diffeomorphisms with that property are called
{\em symplectomorphisms}. Since the volume form of a
symplectic manifold depends only of the canonical $2$-form,
symplectomorphisms preserve volume. We compose with a scaling by 
$\frac{1}{2}$ in the first $n$ variables, that divides $\vol (U)$ by $2^n$, 
and we are done. \qed 

Before proving Proposition~\ref{convexity}, we will need
the following result about convexity which has been attributed to 
Legendre. (See also \cite[Convexity Theorem 1.2]{GROMOV} and a 
generalization in \cite[Th.~5.1]{Avriel}.) 
\begin{leg} 
 If $f$ is convex and of class $\mathcal C^2$ on $\mathbb R^n$, 
then the closure of the image
 $\{ \nabla f _{r} : r \in \mathbb R^n\}$ in $\Rn$ is convex.
\end{leg}

By replacing $f$ by $g_A$, we conclude that the image
of the momentum map $\nabla g_A$ is convex.

\noindent
{\bf Proof of Proposition~\ref{convexity}:} 
The momentum map $\nabla g_A$ maps $\toric{n}$ onto the interior
of $\cvx (A)$. Indeed, let $a=A^{\alpha}$ be a row of $A$, associated to a
vertex of $\cvx (A)$. Then there is a direction $v \in \mathbb R^n$ 
such that
\[
a\cdot v = \max_{x \in \cvx (A)} x \cdot v
\]
for some unique $a$. 

We claim that $a \in \overline{\nabla g_A (\mathbb R^n)}$. Indeed,
let $x(t) = v_A (tv)$, $t$ a real parameter. If $b$ is another 
row of $A$,
\[
e^{a\cdot tv} = e^ {t {a \cdot v}} \gg e^{t{b \cdot v}} = e^{b \cdot tv}
\]
as $t \rightarrow \infty$. We can then write $\hat v_A(tv)^{2T}$ as:
\[
\hat v_A (tv) = 
\left[
\begin{matrix}
\vdots \\
e^{t a \cdot v} \\
\vdots
\end{matrix}
\right]^T
C
\diag
\left[
\begin{matrix}
\vdots \\
e^{t a \cdot v} \\
\vdots
\end{matrix}
\right]
\punct{.}
\]

Since $C$ is positive definite, $C_{\alpha \alpha}>0$ and
\[
\lim_{t \rightarrow \infty}
v_A (tv) ^{2T}
=
\lim_{t \rightarrow \infty}
\frac{ \hat v_A (tv)^{2T} }{\| \hat v_A (tv) \|^2} 
=
\mathrm{e}_a ^T \frac{C_{\alpha \alpha}}{C_{\alpha \alpha}}
=
\mathrm{e}_a ^T
\punct{,}
\]
where $\mathrm{e}_a$ is the unit vector in $\mathbb R^M$ corresponding to
the row $a$. 
It follows that $\lim_{t\rightarrow \infty} \nabla g_A(tv) = a$

When we set $q=0$, we have $\det D^2 g_A \ne 0$ on $\mathbb R^n$, so we
have a local diffeomorphism at each point $p \in \mathbb R^n$.
Assume that $(\nabla g_A)_p = (\nabla g_A)_{p'}$ for $p \ne p'$.
Then, let $\gamma (t) = (1-t)p + tp'$. The function 
$t \mapsto (\nabla g_A)_{\gamma(t)} \gamma'(t)$ has the same value at $0$
and at $1$, hence by Rolle's Theorem its derivative must vanish at some 
$t^* \in (0,1)$. 

In that case,
\[
(D^2 g_A)_{\gamma(t^*)} (\gamma'(t^*), \gamma'(t^*)) = 0
\]
and since $\gamma'(t^*) = p'-p \ne 0$, $\det D^2g_A$ must vanish
in some $p \in \mathbb R^n$. This contradicts
Lemma~\ref{dimnondeg}. \qed  

\subsection{The Condition Matrix}
\label{sec:bernshtein}
Following~\cite{BCSS}, we look at the linearization of the
implicit function $p+q\sqrt{-1} = G(f)$ for the equation
$\mathit{ev}_\fA(f,p+q\sqrt{-1})=0$.

\begin{definition} The {\em condition matrix} of $\mathit{ev}$ at
$(f, p+q\sqrt{-1})$ is
\[
DG = D_{\scriptscriptstyle \mathcal T^n}
(\mathit{ev})^{-1} D_{\scriptscriptstyle \mathcal F}
(\mathit{ev})
\punct{,}
\]
where $\mathcal F = \mathcal F_{A_1} \times \cdots \times \mathcal F_{A_n}$.
\end{definition}

Above, $D_{\scriptscriptstyle \mathcal T^n} (\mathit{ev})$ is a linear operator
from an $n$-dimensional complex space into $\mathbb C^n$, while 
$D_{\scriptscriptstyle \mathcal F} (\mathit{ev})$ goes 
from an $(M_1+\cdots+M_n)$-dimensional complex space into $\mathbb C^n$.

\begin{lemma} \label{herecomesthewedge}
If $p+iq\!\in\cT^n$ and $f(\exp(p+iq))=\bO$ then 
\begin{multline*}
\det \left( DG DG^H \right) ^{-1}
dp_1\wedge dq_1 \wedge \cdots \wedge dp_n\wedge dq_n 
=
(-1)^{n(n-1)/2}\ 
\bigwedge \sqrt{-1} 
f^i \cdot (Dv_{A_i})_{\scriptscriptstyle (p,q)}dp
\wedge \\
\wedge \bar f^i \cdot (Dv_{A_i})_{\scriptscriptstyle (p,-q)}dq
\punct{.}
\end{multline*}
\end{lemma}

\noindent 
Note that although $f^i \cdot (Dv_{A_i})_{\scriptscriptstyle (p,q)}dp$
is a complex-valued form, each wedge
$f^i \cdot (Dv_{A_i})_{\scriptscriptstyle (p,q)}dp
\wedge \bar f^i \cdot (Dv_{A_i})_{\scriptscriptstyle (p,-q)}dq$
is a real-valued $2$-form.

\noindent 
{\bf Proof of Lemma \ref{herecomesthewedge}:} We compute:
\[
D_{\scriptscriptstyle \mathcal F}(\mathit{ev})|_{(p,q)} 
= 
\left[
\begin{matrix}
\sum_{\alpha=1}^{M_1} \hat v_{A_1}^{\alpha} (p+q \sqrt{-1}) df^1_{\alpha}\\
\vdots \\
\sum_{\alpha=1}^{M_n} \hat v_{A_n}^{\alpha} (p+q \sqrt{-1}) df^n_{\alpha}
\end{matrix}
\right]
\punct{,}
\]

and hence
\[
D_{\scriptscriptstyle \mathcal F}
(\mathit{ev}) D_{\scriptscriptstyle \mathcal F}(\mathit{ev})^H = 
\mathrm{diag \ } \| \hat v_ {A_i} \|^2. 
\]

Also,

\[
D_{\scriptscriptstyle \mathcal T^n}(\mathit{ev}) = 
\left[
\begin{matrix}
f^1 \cdot D\hat v_{A_1}\\
\vdots \\
f^n \cdot D\hat v_{A_n}
\end{matrix}. 
\right]
\]

Therefore,

\[
\det \left(DG_{\scriptscriptstyle (p,q)} DG_{\scriptscriptstyle (p,q)}^H\right)^{-1}
=
\left| \det 
\left[
\begin{matrix}
f^1 \cdot \frac{1}{\|\hat v_{A_1}\|} D\hat v_{A_1} \\
\vdots \\
f^n \cdot \frac{1}{\|\hat v_{A_n}\|} D\hat v_{A_n} 
\end{matrix}
\right]
\right|^2. 
\]

We can now use Lemma~\ref{lem-orthog} to conclude the following:

\formula{\label{det:cond:matrix}Determinant of the Condition Matrix}
{
\det \left(DG_{\scriptscriptstyle (p,q)} DG_{\scriptscriptstyle (p,q)}^H\right)^{-1}
=
\left| \det 
\left[
\begin{matrix}
f^1 \cdot Dv_{A_1}\\ 
\vdots \\
f^n \cdot Dv_{A_n} 
\end{matrix}
\right]
\right|^2
}

We can now write the same formula as a determinant of a block matrix:

\[
\det \left(DG_{\scriptscriptstyle (p,q)} DG_{\scriptscriptstyle (p,q)}^H\right)^{-1}
=
\det 
\left[
\begin{matrix}
f^1 \cdot Dv_{A_1}&\\ 
\vdots &\\
f^n \cdot Dv_{A_n}&\\ 
&\bar f^1 \cdot D\bar v_{A_1}\\ 
&\vdots \\
&\bar f^n \cdot D\bar v_{A_n} 
\end{matrix}
\right]
\]

\noindent
and replace the determinant by a wedge.
The factor $(-1)^{n(n-1)/2}$ comes from replacing
$dp_1 \wedge \cdots \wedge dp_n \wedge dq_1 \wedge \cdots \wedge dq_n$
by $dp_1 \wedge dq_1 \wedge \cdots \wedge dp_n \wedge dq_n$. \qed 

We are now ready to prove our main theorems. 

\section{The Proofs of Theorems 1--6}  

We will prove our main theorems in the following order: 
1, 2, 4, 5, 3, 6. 

\subsection{The Proof of Theorem \ref{thm:simpleunmixed}}
The first assertion, modulo an exponential change of coordinates and 
using the multi-projective metric $d_\Pro(\cdot,\cdot)$, follows immediately 
from Theorem \ref{thm:condnumber}. 

As for the rest of Theorem \ref{thm:simpleunmixed}, Theorem \ref{thm:condnumber} applied 
to the {\bf linear} case then provides the following interpretation of $\nu^{\linear} (n, 
\eps)$:
\[
\nu^{\linear} (n, \eps) = 
\mathrm{Prob} \left[ 
d_{\mathbb P}(f, \Sigma_{(p,q)}) < \eps
\right], 
\]
where $f$ is a complex random linear polynomial system,
and $(p,q)$ is such that $f(\exp(p+iq))=0$. So we are on our way to 
proving the inequality \[\mathrm{Prob} \left[
d_{\mathbb P}(f, \Sigma_{(p,q)}) < \eps 
\right]\!\leq\!n^3(n+1)\vol(A)(\#A-1)(\#A-2)\eps^4,\] 
for {\bf general} $f$, which clearly implies our desired bound. 

To prove the latter inequality, recall that by the definition of the multi-projective distance 
$d_{\mathbb P}(\cdot,\cdot)$, we have the following equality: 
\[
d_{\mathbb P}(f, \Sigma_{(p,q)})^2
=
\min_{\substack{g \in \Sigma_{(p,q)}\\ \lambda \in (\mathbb C^*)^n}}
\sum_{i=1}^n \frac{\|f^i - \lambda_i g^i\|^2}{\|f^i\|^2}. 
\]
So let $g$ be so that the above minimum is attained. Without loss of
generality, we may scale the $g^i$ so that $\lambda_1 =
\cdots = \lambda_n = 1$. In that case,
\[
d_{\mathbb P}(f, \Sigma_{(p,q)})^2
=
\sum_{i=1}^n \frac{\|f^i - g^i\|^2}{\|f^i\|^2}
\ge
\frac{ \sum_{i=1}^n \|f^i - g^i\|^2}{\sum_{j=1}^n \|f^j\|^2}. 
\]

\noindent
We are then in the setting of \cite[pp.\ 248--250]{BCSS} where we 
identify our linear $f$ with a normally distributed $(n+1)\times n$ complex matrix. 
The right-hand side in the above inequality is then precisely the left-hand term in 
\cite[Rem.\ 2, Pg.\ 250]{BCSS}.  Therefore, using the notation of \cite[Prop.\ 4]{BCSS}, 
$d_{\mathbb P}(f, \Sigma_{(p,q)}) \ge d_{\mathrm{F}}(f, \Sigma_x)$. So 
it follows that
\[
\nu^{\linear} (n, \eps) = 
\mathrm{Prob} \left[ 
d_{\mathbb P}(f, \Sigma_{(p,q)}) < \eps
\right]
\le
\mathrm{Prob} \left[ 
d_{\mathrm F}(f, \Sigma_{x}) < \eps
\right]
\]
\noindent 
and the last probability is bounded above by 
$n^3 (n+1)(\#A-1)(\#A-2)\eps^4$ 
via \cite[Thm.\ 6, Pg.\ 254]{BCSS}. Theorem \ref{thm:simpleunmixed} now 
follows. \qed

\subsection{The Proof of Theorem \ref{GEN-BERNSHTEIN}} 
Using \cite[Theorem 5 p.~243]{BCSS} (or Proposition~\ref{coarea2},  
Pg.\ \pageref{coarea2} below), we deduce that the average number of 
complex roots is:
\[
\mathrm{Avg} =
\int_{(p,q) \in U}
\int_{f \in \mathcal F_{(p,q)}}
\left(\prod \frac{e^{-\|f^i\|^2/2}}{(2 \pi)^{M_i}}\right) 
\det \left(DG_{\scriptscriptstyle (p,q)} DG_{\scriptscriptstyle (p,q)}^H\right)^{-1}. 
\]

By Lemma~\ref{herecomesthewedge}, we can replace the inner integral by
a $2n$-form valued integral:
\begin{multline*}
\mathrm{Avg} =
(-1)^{n(n-1)/2}\ 
\int_{(p,q) \in U}
\int_{f \in \mathcal F_{(p,q)}}
\bigwedge_i \frac{e^{-\|f^i\|^2/2}}{(2 \pi)^{M_i}}
f^i \cdot (Dv_{A_i})_{\scriptscriptstyle (p,q)}dp
\wedge \\
\wedge \bar f^i \cdot (Dv_{A_i})_{\scriptscriptstyle (p,-q)}dq
\punct{.}
\end{multline*}

Since the image of $Dv_{A_i}$ is precisely $\mathcal (F_{A_i})_{(p,q)} 
\subset \mathcal F_{A_i}$,
one can add $n$ extra variables corresponding to the directions 
$v_{A_i}(p+q \sqrt{-1})$
without changing the integral: we write $\mathcal F_{A_i} = \mathcal 
F_{A_i,(p,q)} \times
\mathbb C v_{A_i}(p+q \sqrt{-1})$. Since $\left( f^i + t v_{A_i}(p+q \sqrt{-1}) \right)
Dv_{A_i}$ is equal to $f^i Dv_{A_i}$, the average number of roots is indeed:
\begin{multline*}
\mathrm{Avg} =
(-1)^{n(n-1)/2}\ 
\int_{(p,q) \in U}
\int_{f \in \mathcal F}
\bigwedge_i \frac{e^{-\|f^i\|^2/2}}{(2 \pi)^{M_i +1}}
f^i \cdot (Dv_{A_i})_{\scriptscriptstyle (p,q)}dp
\wedge \\
\wedge \bar f^i \cdot (Dv_{A_i})_{\scriptscriptstyle (p,-q)}dq
\punct{.}
\end{multline*}

In the integral above, all the terms that are multiple of $f^i_{\alpha} \bar 
f^i_{\beta}$
for some $\alpha \ne \beta$ will cancel out. Therefore,
\begin{multline*}
\mathrm{Avg} =
(-1)^{n(n-1)/2}\ 
\int_{(p,q) \in U}
\int_{f \in \mathcal F}
\bigwedge_i \frac{e^{-\|f^i\|^2/2}}{(2 \pi)^{M_i+1}}
\sum_{\alpha} |f^i_{\alpha}|^2 (Dv_{A_i})^{\alpha}_{\scriptscriptstyle (p,q)}dp
\wedge 
\\
\wedge
(Dv_{A_i})^{\alpha}_{\scriptscriptstyle (p,-q)}dq
\punct{.}
\end{multline*}

Now, we apply the integral formula:
\[
\int_{x \in \mathbb C^{M} } |x_1|^2 \frac{e^{-\|x\|^2 / 2}}{(2 \pi)^{M}}
=
\int_{x_1 \in \mathbb C} |x_1|^2 \frac{e^{-|x_1|^2 / 2}}{2 \pi}
= 2
\]
to obtain:
\[
\mathrm{Avg} =
\frac{(-1)^{n(n-1)/2}}{\pi^n}
\int_{(p,q) \in U}
\bigwedge \sum_{\alpha} (Dv_{A_i})^{\alpha}_{\scriptscriptstyle (p,q)}dp
\wedge (Dv_{A_i})^{\alpha}_{\scriptscriptstyle (p,-q)}dq
\punct{.}
\]

According to formul\ae~\ref{secder} and~\ref{omega}, the integrand is just 
$2^{-n} \bigwedge \omega_{A_i}$, and thus
\[
\mathrm{Avg} =
\frac{(-1)^{n(n-1)/2}}{\pi^n}
\int_{U}
\bigwedge_i \omega_{A_i}
=
\frac{n!}{\pi^n}
\int_{U}
d\toric{n}. \text{\qed}
\]

\subsection{The Proof of Theorem \ref{thm:condnumber}}

Let $(p,q) \in \toric{n}$ and let
$f \in \mathcal F_{(p,q)}$. Without loss of generality, we can
assume that $f$ is scaled so that for all $i$,
$\|f^i\|=1$. 

Let $\delta f \in \mathcal F_{(p,q)}$
be such that $f+\delta f$ is singular at $(p,q)$, and assume
that $\sum \| \delta f^i \|^2$ is minimal. Then, due
to the scaling we choose,
\[
d_{\mathbb P} (f, \Sigma_{(p,q)}) = \sqrt{\sum \| \delta f^i \|^2} 
\punct{.}
\]

Since $f + \delta f$ is singular, there is a vector $u \ne 0$ such
that
\[
\left[
\begin{matrix}
(f^1 + \delta f^1) \cdot (D\hat v_{A_1})_{(p,q)}
\\
\vdots
\\
(f^n +\delta f^n) \cdot (D\hat v_{A_n})_{(p,q)}
\end{matrix}
\right]
u = 0
\]
and hence
\[
\left[
\begin{matrix}
(f^1 + \delta f^1) \cdot (Dv_{A_1})_{(p,q)}
\\
\vdots
\\
(f^n +\delta f^n) \cdot (Dv_{A_n})_{(p,q)}
\end{matrix}
\right]
u = 0
\punct{.}
\]
\medskip
 
This means that
\[
\left\{
\begin{array}{lcl}
f^1 \cdot Dv_{A_1} u &=& - \delta f^1 \cdot Dv_{A_1} u \\
&\vdots&\\
f^n \cdot Dv_{A_n} u &=& - \delta f^n \cdot Dv_{A_n} u \\
\end{array}
\right.
\punct{.}
\]

Let $D(f)$ denote the matrix
\[
D(f) \defeq
\left[
\begin{matrix}
f^1 \cdot (Dv_{A_1})_{(p,q)}
\\
\vdots
\\
f^n \cdot (Dv_{A_n})_{(p,q)}
\end{matrix}
\right]
\punct{.}
\]

Given $v = D(f) \ u$, we obtain:
\begin{equation}\label{eqonv}
\left\{
\begin{array}{lcl}
v_1 &=& - \delta f^1 \cdot Dv_{A_1} D(f)^{-1} v \\
&\vdots&\\
v_n &=& - \delta f^n \cdot Dv_{A_n} D(f)^{-1} v \\
\end{array}
\right.
\end{equation}

We can then scale $u$ and $v$, such that $\|v\|=1$.

\begin{claim}
Under the assumptions above, $\delta f^i$ is colinear to 
$\left( Dv_{A_i} D(f)^{-1} v \right)^H$.
\end{claim}

\noindent {\bf Proof:}
Assume that $\delta f^i = g + h$, with $g$ colinear and $h$ orthogonal to
$\left( Dv_{A_i} D(f)^{-1} v \right)^H$.
As the image of $Dv_{A_i}$ is orthogonal to $v_{A_i}$, $g$ is orthogonal
to $v_{A_i}^H$, so $\mathit{ev}(g^i,(p,q))=0$ and hence
$\mathit{ev}(h^i,(p,q))=0$. We can therefore
replace $\delta f^i$ by $g$ without compromising
equality~(\ref{eqonv}). Since $\|\delta f\|$ was minimal, this implies
$h=0$.
\qed

We obtain now an explicit expression for $\delta f^i$ in terms of $v$:

\[
\delta f^i
= - v_i
\frac{
 \left( Dv_{A_i} D(f)^{-1} v \right)^H
}
{
\| Dv_{A_i} D(f)^{-1} v \|^2
}
\punct{.}
\]

Therefore,
\[
\| \delta f^i \|
=
\frac{ |v_i| }
{\|Dv_{A_i} D(f)^{-1} v \|}
=
\frac{ |v_i| }
{\|\left( D(f)^{-1} v \right)\|_{A_i}}
\punct{.}
\]

So we have proved the following result: 

\begin{lemma} \label{lemonv}
Fix $v$ so that $\|v\|=1$ and 
let $\delta f \in \mathcal F_{(p,q)}$ be such that equation (\ref{eqonv}) 
holds and $\|\delta f\|$ is minimal. Then, 
\[
\| \delta f^i \|
=
\frac{ |v_i| }
{\|D(f)^{-1} v \|_{A_i}}
\punct{.}
\]
\end{lemma}

Lemma~\ref{lemonv} provides an immediate lower bound for 
$\| \delta f \| = \sqrt{\sum \| \delta f^i\|^2}$: Since
\[
\| \delta f^i \| \ge
\frac{ |v_i| }
{\max_j \| D(f)^{-1} v \|_{A_j}}
\punct{,}
\]
we can use $\|v\|=1$ to deduce that
\[
\sqrt{\sum_i \| \delta f^i \|^2} \ge
\frac{ 1 }
{\max_j \|D(f)^{-1} v \|_{A_j}}
\ge
\frac{ 1 }
{\max_j \| D(f)^{-1} \|_{A_j}}
\punct{.}
\]

Also, for any $v$ with $\|v\|=1$, we can choose $\delta f$ minimal so
that equation~(\ref{eqonv}) applies. Using Lemma~\ref{lemonv}, we
obtain:
\[
\| \delta f^i \| \le
\frac{ |v_i| }
{\min_j \| D(f)^{-1} v \|_{A_j}}
\punct{.}
\]

Hence
\[
\sqrt{\sum_i \| \delta f^i \|^2} \le
\frac{ 1 }
{\min_j \| D(f)^{-1} v \|_{A_j}}
\punct{.}
\]

Since this is true for any $v$, and $\| \delta f\|$ is
minimal for all $v$, we have 
\[
\sqrt{\sum_i \| \delta f^i \|^2} \le
\frac{ 1 }
{\max_{\|v\|=1} \min_j \| D(f)^{-1}\|_{A_j}}
\]

and this proves Theorem \ref{thm:condnumber}.

\subsection{The Idea Behind the Proof of Theorem~\ref{mixed-cond}}

The proof of Theorem~\ref{mixed-cond} is long. 
We first sketch the idea of the proof. Recall
that $\mathcal F_{(p,q)}$ is the set of all
$f \in \mathcal F$ such that $\mathit{ev}(f;p
+q\sqrt{-1}) = 0$, and that $\Sigma_{(p,q)}$
is the restriction of the discriminant to the
fiber $\mathcal F_{(p,q)}$:
\[
\Sigma_{(p,q)} \defeq \{ f \in \mathcal F_{(p,q)}:
 D(f)_{(p,q)} \text{\ does not have full rank} \} 
\punct{.}
\]
\medskip

The space $\mathcal F$ is endowed with a Gaussian
probability measure, with volume element 
\[
\frac{e^{-\|f\|^2 / 2}}{(2 \pi)^{\sum M_i}} d\mathcal F
\punct{,}
\]
where $d\mathcal F$ is the usual volume form in 
$\mathcal F = (\mathcal F_{A_1}, \langle\cdot,\cdot
\rangle_{\scriptscriptstyle A_1}) \times \cdots
\times (\mathcal F_{A_n}, \langle\cdot,\cdot\rangle_{\scriptscriptstyle A_n})$
and $\|f\|^2 = \sum \|f^i\|_{A_i}^2$.
For $U$ a set in $\toric{n}$, we defined earlier (in the statement of 
Theorem \ref{mixed-cond}) the quantity: 
\[
\nu^{A} (U,\eps) \defeq 
\mathrm{Prob} [ \mymu(f,U) > \eps^{-1}]
=
\mathrm{Prob} [ 
\exists (p,q) \in U: d_{\mathbb P}(f,\Sigma_{(p,q)}) < \eps]
\punct{.}
\]

The na\"{\i}ve idea for bounding $\nu^A(U,\eps)$ is as follows:
Let $V(\eps) \defeq \{ (f,(p,q)) \in \mathcal F \times U: 
\mathit{ev}(f;(p,q))=0 \text{\ and \ }
d_{\mathbb P} (f, \Sigma_{(p,q)}) < \eps \}$. We also define $\pi: V(\eps)
\rightarrow \mathcal F$ as the canonical projection mapping 
$\cF\times U$ to $\cF$, and set 
$\#_{V(\eps)}(f) \defeq \# \{ (p,q) \in U: (f,(p,q)) \in V(\eps) \}$.
Then,

\begin{eqnarray*}
\nu^A(U,\eps)
&=&
\int_{f \in \mathcal F}
\chi_{\pi(V(\eps))} (f) \ 
\frac{e^{-\|f\|^2 / 2}}{(2 \pi)^{\sum M_i}} d\mathcal F
\\
&\le&
\int_{f \in \mathcal F}
\#_{V(\eps)} \ 
\frac{e^{-\|f\|^2 / 2}}{(2 \pi)^{\sum M_i}} d\mathcal F
\end{eqnarray*}
with equality in the linear case.

Now we apply the coarea formula~\cite[Theorem 5 p.~243]{BCSS} to obtain:

\[
\nu^A(U,\eps)
\le
\int_{(p,q) \in U \subset \toric{n}}
\int_{
\substack{
f \in \mathcal F_{(p,q)} \\
d_{\mathbb P} (f, \Sigma_{(p,q)}) < \eps
}
}
\frac{1}{NJ(f;(p,q))} \ 
\frac{e^{-\|f\|^2 / 2}}{(2 \pi)^{\sum M_i}} d\mathcal F
\ dV_{\toric{n}}
\punct{,}
\]

where $dV_{\toric{n}}$ stands for Lebesgue measure in
$\toric{n}$. Again, in the linear case, we have
equality.

We already know from Lemma~\ref{herecomesthewedge}
that 
\[
1/NJ(f;(p,q))
=
\bigwedge_{i=1}^n 
f^i \cdot (Dv_{A_i})_{(p,q)} dp \wedge \bar f^i \cdot (D\bar v_{A_i})_{(p,q)} dq
\punct{.}
\]

We should focus now on the inner integral. In each coordinate space
$\mathcal F_{A_i}$, we can introduce a new orthonormal system of coordinates
(depending on $(p,q)$) by decomposing:
\[
f^i = f^i_{\romone} + f^i_{\romtwo} + f^i_{\romthree}
\punct{,}
\]

where $f^i_{\romone}$ is the component colinear
to $v_{A_i}^H$, 
$f^i_{\romtwo}$ is the projection of
$f^i$ to $({\mathrm{range}}\  Dv_{A_i})^H$,
and
$f^i_{\romthree}$ is orthogonal to 
$f^i_{\romone}$ and $f^i_{\romtwo}$.

Of course, $f^i \in (\mathcal F_{A_i})_{(p,q)}$ if and only if
$f^i_{\romone} = 0$.

Also,
\begin{multline*}
\bigwedge_{i=1}^n
f^i \cdot (Dv_{A_i})_{(p,q)} dp \wedge \bar f^i \cdot (D\bar v_{A_i})_{(p,q)} dq
=
\\
=
\bigwedge_{i=1}^n
f_{\romtwo}^i \cdot (Dv_{A_i})_{(p,q)} dp \wedge \bar f_{\romtwo}^i \cdot 
(D\bar v_{A_i})_{(p,q)} dq \punct{.}
\end{multline*}

It is an elementary fact that 
\[
d_{\mathbb P}(f^i_{\romtwo}+f^i_{\romthree}, \Sigma_{(p,q)})
\le
d_{\mathbb P}(f^i_{\romtwo}, \Sigma_{(p,q)})
\punct{.}
\]

It follows that for $f \in \mathcal F_{(p,q)}$:
\[
d_{\mathbb P}(f, \Sigma_{(p,q)})
\le
d_{\mathbb P}(f_{\romtwo}, \Sigma_{(p,q)})
\punct{,}
\]
with equality in the linear case. Hence, we obtain:

\begin{multline*}
\nu^A(U,\eps)
\le
\int_{(p,q) \in U \subset \toric{n}}
\int_{
\substack{
f \in \mathcal F_{(p,q)} \\
d_{\mathbb P} (f_{\romtwo}, \Sigma_{(p,q)}) < \eps
}}
\left(
\bigwedge_{i=1}^n
f_{\romtwo}^i \cdot (Dv_{A_i})_{(p,q)} dp \wedge \bar f_{\romtwo}^i \cdot (D\bar v_{A_i})_{(p,q)} dq
\right)
\cdot
\\
\cdot
\frac{e^{-\|f^i_{\romtwo}+f^i_{\romthree}\|^2 / 2}}{(2 \pi)^{\sum M_i}} d\mathcal F
\ dV_{\toric{n}}
\punct{,}
\end{multline*}
with equality in the linear case. We can integrate the 
$\sum (M_i - n - 1)$ variables $f_{\romthree}$
to obtain:

\begin{proposition}\label{assimpleasitgets}
\begin{multline*}
\nu^A(U,\eps)
\le
\int_{(p,q) \in U \subset \toric{n}}
\int_{
\substack
{
f_{\romtwo} \in \mathbb C^{n^2} \\
d_{\mathbb P} (f_{\romtwo}, \Sigma_{(p,q)}) < \eps
}}
\left(
\bigwedge_{i=1}^n
f_{\romtwo}^i \cdot (Dv_{A_i})_{(p,q)} dp \wedge \bar f_{\romtwo}^i \cdot (D\bar v_{A_i})_{(p,q)} dq
\right)
\cdot \\ \cdot
\frac{e^{-\|f^i_{\romtwo}\|^2 / 2}}{(2 \pi)^{n(n+1)}} 
\ dV_{\toric{n}}
\punct{.}
\end{multline*}
with equality in the linear case. \qed 
\end{proposition}

\subsection{The Proof of Theorem \ref{mixed-cond}}

The domain of integration in Proposition~\ref{assimpleasitgets}
makes integration extremely difficult. In order to estimate
the inner integral, we will need to perform a change of
coordinates.

Unfortunately, the Gaussian in Proposition~\ref{assimpleasitgets}
makes that change of coordinates extremely hard, and we will have
to restate Proposition \ref{assimpleasitgets} in terms of integrals over 
a product of projective spaces.  

The domain of integration will be $\mathbb P^{n-1} \times 
\cdots \times \mathbb P^{n-1}$. Translating an integral
in terms of Gaussians to an integral in terms of
projective spaces is not immediate, and we will use
the following elementary fact about Gaussians:

\begin{lemma}\label{gaussians1}
 Let $\varphi: \mathbb C^n \rightarrow \mathbb R$ be
$\mathbb C^*$-invariant (in the sense of the usual 
scaling action). Then we can also interpret
$\varphi$ as a function from $\mathbb P^{n-1}$ into
$\mathbb R$, and:
\[
\frac{1}{\vol(\mathbb P^{n+1})} 
\int_{[x]\in \mathbb P^{n-1}} \varphi(x)d[x]  =
\int_{x \in \mathbb C^{n}}  \varphi(x) 
\frac{e^{-\|x\|^2/2}}{(2\pi)^n} dx 
\punct{,}
\]
where, respectively, the natural volume forms on $\Pro^{n-1}$ and $\Cn$ 
are understood for each integral. 
\end{lemma}

Now the integrand in Proposition~\ref{assimpleasitgets}
is not $\mathbb C^*$-invariant. This is why we will need 
the following formula: 

\begin{lemma}\label{gaussians2}
 Under the hypotheses of Lemma~\ref{gaussians1},
\[
\frac{1}{\vol(\mathbb P^{n+1})} 
\int_{[x]\in \mathbb P^{n-1}} \varphi(x) d[x] =
\frac{1}{2n}
\int_{x \in \mathbb C^{n}}  \|x\|^2 \varphi(x) 
\frac{e^{-\|x\|^2/2}}{(2\pi)^n} dx
\punct{.}
\]
where, respectively, the natural volume forms on $\Pro^{n-1}$ and $\Cn$ 
are understood for each integral. 
\end{lemma}

\noindent {\bf Proof:}
\begin{eqnarray*}
\int_{x \in \mathbb C^{n}}  \|x\|^2 \varphi(x) 
\frac{e^{-\|x\|^2/2}}{(2\pi)^n} dx  
&=&
\int_{\Theta \in S^{2n-1}} 
\int_{r=0}^{\infty}
|r|^{2n+1} \varphi(\Theta) 
\frac{e^{-|r|^2/2}}{(2\pi)^n} dr d\Theta 
\end{eqnarray*} 
\[=
\int_{\Theta \in S^{2n-1}} 
\left(
-
\left[
|r|^{2n} 
\frac{e^{-|r|^2/2}}{(2\pi)^n} 
\right]_0^{\infty}
\right. 
\left. 
+
2n
\int_{r=0}^{\infty}
|r|^{2n-1} 
\frac{e^{-|r|^2/2}}{(2\pi)^n} 
dr 
\right)
\varphi(\Theta) 
d\Theta \] 
\[=
2n
\int_{x \in \mathbb C^{n}} \varphi(x) 
\frac{e^{-\|x\|^2/2}}{(2\pi)^n} dx 
\] 
\qed

We can now introduce the notation:

\[
\mathrm{WEDGE}^A(f_{\romtwo}) 
\defeq 
\bigwedge_{i=1}^n 
\frac{1}{\|f_{\romtwo}^i\|^2} 
f_{\romtwo}^i \cdot (Dv_{A_i})_{(p,q)} dp \wedge 
\bar f_{\romtwo}^i \cdot (D\bar v_{A_i})_{(p,q)} dq
\punct{.}
\]

This function is invariant under the $(\mathbb C^*)^n$-action
$\lambda \star f_{\romtwo} : f_{\romtwo} \mapsto 
(\lambda_1 f_{\romtwo}^1, \cdots, \lambda_n f_{\romtwo}^n)$. 

We adopt the following conventions:
$\mathcal F_{\romtwo} \subset \mathcal F$ is the space spanned by
coordinates $f_{\romtwo}$ and $\mathbb P (\mathcal F_{\romtwo}) $
is its quotient by $(\mathbb C^*)^n$.

We apply $n$ times
Lemma~\ref{gaussians2} and obtain:

\begin{proposition}\label{multiprojective}
Let $\mathrm{VOL} \defeq \vol(\mathbb P^{n-1})^{n}$. Then,
\[
\nu^A(U,\eps)
\le
\frac{(2n)^n}{\mathrm{VOL}}
\int_{(p,q) \in U \subset \toric{n}}
\int_{
\substack{
f_{\romtwo} \in \mathbb P (\mathcal F_{\romtwo})\\
d_{\mathbb P} (f_{\romtwo}, \Sigma_{(p,q)}) < \eps
} }
\mathrm{WEDGE}^A(f_{\romtwo}) \
d\mathbb P (\mathcal F_{\romtwo})
\ dV_{\toric{n}}
\]
and in the linear case,
\[
\nu^{\linear}(U,\eps)
=
\frac{(2n)^n}{\mathrm{VOL}}
\int_{(p,q) \in U \subset \toric{n}}
\int_{
\substack{
g_{\romtwo} \in \mathbb P (\mathcal F_{\romtwo}^{\linear}) \\
d_{\mathbb P} (g_{\romtwo}, \Sigma_{(p,q)}^{\linear}) < \eps
} }
\mathrm{WEDGE}^{\linear}(g_{\romtwo}) \
d\mathbb (P \mathcal F_{\romtwo}^{\linear})
dV_{\toric{n}} \ \text{\qed} 
\]
\end{proposition}

Now we introduce the following change of coordinates. Let
$L \in GL(n)$ be such that the minimum in Definition~\ref{intrinsic}
p.~\pageref{intrinsic} is attained:
\[
\begin{array}{crcl}
 \varphi: &
\mathbb P^{n-1} \times \cdots \times \mathbb P^{n-1} &
\rightarrow &
\mathbb P^{n-1} \times \cdots \times \mathbb P^{n-1} \\
&
f_{\romtwo} & \mapsto & g_{\romtwo} \defeq \varphi(f_{\romtwo}) \punct{,}\text{\ such}\\
& & & \text{that\ } 
g_{\romtwo}^i = f_{\romtwo}^i \cdot Dv_{A_i} L \punct{.}
\end{array}
\]

Without loss of generality, we scale $L$ such that $\det L=1$.
The following property follows from the definition of
$\mathrm{WEDGE}$:

\begin{equation}\label{wedge1}
\mathrm{WEDGE}^A (f_{\romtwo})
=
\mathrm{WEDGE}^{\linear} (g_{\romtwo}) \ 
\prod_{i=1}^n
\frac {\| g_{\romtwo}^i \|^2} {\| f_{\romtwo}^i \|^2}
\end{equation}

Assume now that
$d_{\mathbb P} (f_{\romtwo}, \Sigma_{(p,q)}) < \eps$.
Then there is $\delta f \in \mathcal F_{\romtwo}$, 
such that $f + \delta f \in \Sigma_{(p,q)}^{\linear}$ and
$\| \delta f \| \le \eps$ (assuming the scaling  $\|f^i_{\romtwo}\|=1$
for all $i$).

Setting $g_{\romtwo} = \varphi (f_{\romtwo})$ and $\delta g = \varphi(g)$,
we obtain that $g + \delta g \in \Sigma_{(p,q)}^{\linear}$. 
\[
d_{\mathbb P} (g, \Sigma_{(p,q)}^{\linear})
\le
\sqrt{ \sum_{i=1}^n \frac {\| \delta g^i \|^2}{\| g^i_{\romtwo} \|^2}
}
\]

At each value of $i$,
\[
\frac {\| \delta g^i \|}{\| g^i_{\romtwo} \|}
\le
\frac {\| \delta f^i \|}{\| f^i_{\romtwo} \|}
\kappa ( D_{f^i_{\romtwo}}\varphi^i) 
\]
where $\kappa$ denotes Wilkinson's condition number of the
linear operator $D_{f^i_{\romtwo}}\varphi^i$. This is precisely
$\kappa (Dv_{A_i} L)$. Thus,
\[
d_{\mathbb P} (g, \Sigma_{(p,q)}^{\linear})
\le
\eps
\max_i \kappa (Dv_{A_i} L) = \max_i \sqrt{\kappa ( \omega_{A_i} )}
\]
 
Thus, an $\eps$-neighborhood of $\Sigma^A_{(p,q)}$ is mapped
into a $\sqrt{\kappa_U} \eps$ neighborhood of
$\Sigma^{\linear}_{(p,q)}$.

We use this property and equation~(\ref{wedge1}) to bound:

\begin{multline}\label{wedge2}
\nu^A(U,\eps)
\le
\frac{(2n)^n}{\mathrm{VOL}}
\int_{(p,q) \in U \subset \toric{n}}
\int_{
\substack{
g_{\romtwo} \in \mathbb P^{n-1} \times \cdots \times \mathbb P^{n-1} \\
d_{\mathbb P} (g_{\romtwo}, \Sigma^{\linear}_{(p,q)}) < \sqrt{\kappa_U} \eps
} }
\mathrm{WEDGE}^{\linear}(g_{\romtwo})
\cdot
\\
\cdot
\prod_{i=1}^n
\frac {\| g_{\romtwo}^i \|^2} {\| f_{\romtwo}^i \|^2}
|J_{g_{\romtwo}} \varphi^{-1}|^2
\
d(\mathbb P^{n-1}\times \cdots \times \mathbb P^{n-1})
\ dV_{\toric{n}}
\end{multline}

where $J_{g_{\romtwo}} \varphi^{-1}$ is the Jacobian of $\varphi^{-1}$ at
$g_{\romtwo}$. 
\begin{rem}
Considering each $Dv_{A_i}$ as a map from $\mathbb C^n$
into $\mathbb C^n$, the Jacobian is:
\[
J_{g_{\romtwo}} \varphi^{-1}
=
\prod_{i=1}^n \frac{\|\varphi^{-1}(g_{\romtwo})^i\|^n}{\|g_{\romtwo}^i\|^n} 
\left( \det Dv_{A_i}^H Dv_{A_i} \right)^{-1/2}
\punct{.}
\]
We will not use this value in the sequel. \dia  
\end{rem}

In order to simplify the expressions for the bound on $\nu^A (U,\eps)$, it
is convenient to introduce the following notation:

\begin{eqnarray*}
dP
&\defeq&
\frac{(2n)^n}{\mathrm {VOL}}
\mathrm{WEDGE}^{\linear}(g_{\romtwo})
\frac{ d (\mathbb P^{n-1}\times \cdots \times \mathbb P^{n-1}) }
{n! \ (\omega_{\linear})^{\bigwedge n}}
\\
H &\defeq&
\prod_{i=1}^n
\frac {\| g_{\romtwo}^i \|^2} {\| f_{\romtwo}^i \|^2}
|J_g \varphi^{-1}|^2
\\
\chi_{\delta} &\defeq& \chi_{ \left\{ g: d_{\mathbb P}(g, \Sigma_{(p,q)}^{\linear}) < \delta \right\} }
\end{eqnarray*}

Now equation~(\ref{wedge2}) becomes:

\begin{equation}\label{wedge3}
\nu^A(U,\eps)
\le
n!
\int_{(p,q) \in U \subset \toric{n}}
(\omega_{\linear})^{\bigwedge n}
\int_{ g_{\romtwo} \in \mathbb P^{n-1} \times \cdots \times \mathbb P^{n-1} }
dP \
H(g_{\romtwo})
\
\chi_{\sqrt{\kappa_U} \eps} (g_{\romtwo})
\end{equation}

\begin{lemma}
 Let $(p,q)$ be fixed. Then $\mathbb P^{n-1} \times \cdots \times \mathbb P^{n-1}$
together with density  function $dP$, is a probability space.
\end{lemma}

\noindent {\bf Proof:}
  The expected number of roots in $U$ for a linear system is 
\[
n!
\int_{(p,q) \in U}
\omega_{\linear}^{\bigwedge n}
\int_{ g_{\romtwo} \in \mathbb P^{n-1} \times \cdots \times \mathbb P^{n-1} }
dP
\]
which is also $n! \int_U \omega_{\linear}^{\bigwedge n}$. This holds for all $U$,
hence the volume forms are the same and
\[
\int_{ g_{\romtwo} \in \mathbb P^{n-1} \times \cdots \times \mathbb P^{n-1} }
dP
= 1. \text{\qed} \]

This allows us to interpret the inner integral of equation~(\ref{wedge3})
as the expected value of a product. This is less than the product of
the expected values, and:

\begin{multline*}
\nu^A(U,\eps)
\le
n!
\int_{(p,q) \in U \subset \toric{n}}
(\omega_{\linear})^{\bigwedge n}
\left(
\int_{ g_{\romtwo} \in \mathbb P^{n-1} \times \cdots \times \mathbb P^{n-1} }
dP \
H(g_{\romtwo})
\right)
\cdot
\\
\cdot
\left(
\int_{ g_{\romtwo} \in \mathbb P^{n-1} \times \cdots \times \mathbb P^{n-1} }
dP \
\chi_{\sqrt{\kappa_U} \eps} (g_{\romtwo})
\right)
\end{multline*}

Because generic (square) systems of linear equations have exactly one root, we can also
consider $U$ as a probability space, with probability measure
$\frac{1}{\vol^{\linear} (U)} n! \omega_{\linear}^{\bigwedge n}$.
Therefore, we can bound:
\begin{multline*}
\nu^A(U,\eps)
\le
\frac{1}{\vol^{\linear} (U)}
\left(
\int_{(p,q) \in U}
n!
(\omega_{\linear})^{\bigwedge n}
\int_{ g_{\romtwo} \in \mathbb P^{n-1} \times \cdots \times \mathbb P^{n-1} }
dP \
H(g_{\romtwo})
\right)
\cdot \\ \cdot
\left(
\int_{(p,q) \in U}
n!
(\omega_{\linear})^{\bigwedge n}
\int_{ g_{\romtwo} \in \mathbb P^{n-1} \times \cdots \times \mathbb P^{n-1} }
dP \
\chi_{\sqrt{\kappa_U} \eps} (g_{\romtwo})
\right)
\end{multline*}

The first parenthetical expression is $\vol^A (U)$, the volume of $U$
with respect to the toric volume form associated to $A=(A_1, \cdots, A_n)$. 
The second parenthetical expression is $\nu^{\linear}(\sqrt{\kappa_U}\eps,U)$.
This concludes the proof of Theorem~\ref{mixed-cond}. \qed

\subsection{The Proof of Theorem~\ref{expected-real}}
As in the complex case (Theorem~\ref{GEN-BERNSHTEIN}),
the expected number of roots can be computed by
applying the co-area formula:
\[
AVG =
\int_{p \in U}
\int_{f \in \mathcal F_p^{\mathbb R}}
\prod_{i=1}^n 
\frac{e ^{-\|f^i\|^2/2}}{\sqrt{2 \pi}^{M_i}}
\det (DG \ DG^H)^{-1/2}
\punct{.}
\]

Now there are three big diferences. The set $U$ is
in $\mathbb R^n$ instead of $\toric{n}$, the
space $\mathcal F_p^{\mathbb R}$ contains
only real polynomials (and therefore has half the
dimension), and we are integrating the square root
of $1/ \det (DG\ DG^H)$.

Since we do not know in general how to integrate
such a square root, we bound the inner integral
as follows. We consider the real Hilbert space of
functions integrable in $\mathcal F_p^{\mathbb R}$
endowed with Gaussian probability measure. The
inner product in this space is:

\[
\langle \varphi , \psi\rangle \defeq
\int_{\mathcal F_p^{\mathbb R}}
\varphi(f) \psi(f)
\prod_{i=1}^n 
\frac{e^{-\|f^i\|^2/2}}{\sqrt{2 \pi}^{M_i-1}}
dV
\punct{,}
\]

where $dV$ is Lebesgue volume. If $\mathbf 1$ denotes
the constant function equal to $1$, we interpret

\[
AVG
=
\int_{p \in U}
(2\pi)^{-n/2}
\left\langle
\det (DG \ DG^H)^{-1/2}
,
\mathbf 1
\right\rangle
\punct{.}
\]

Hence Cauchy-Schwartz inequality implies:

\[
AVG
\le
\int_{p \in U}
(2\pi)^{-n/2}
\|
\det (DG \ DG^H)^{-1/2}
\|
\|
\mathbf 1
\|
\punct{.}
\]

By construction, $\| \mathbf 1\|=1$, and we are left with:

\[
AVG
\le
\int_{p \in U}
(2\pi)^{-n/2}
\sqrt{
\int_{\mathcal F_p^{\mathbb R}}
\prod_{i=1}^n 
\frac{e^{-\|f^i\|^2/2}}{\sqrt{2 \pi}^{M_i-1}}
\det (DG \ DG^H)^{-1}
}
\punct{.}
\]

As in the complex case, we add extra $n$ variables:
\[
AVG
\le
(2\pi)^{-n/2}
\int_{p \in U}
\sqrt{
\int_{\mathcal F^{\mathbb R}}
\prod_{i=1}^n 
\frac{e^{-\|f^i\|^2/2}}{\sqrt{2 \pi}^{M_i}}
\det (DG \ DG^H)^{-1}
}
\punct{,}
\]

and we interpret $\det (DG \ DG^H)^{-1}$ in terms of a wedge.
Since
\[
\int_{x \in \mathbb R^M} 
|x_1|^2 \frac{e^{-\|x\|^2/2}}{\sqrt{2 \pi}^M}
=
\int_{y \in \mathbb R} 
y^2 \frac{e^{-y^2/2}}{\sqrt{2 \pi}}
=
\int_{y \in \mathbb R} 
\frac{e^{-y^2/2}}{\sqrt{2 \pi}}
=1
\punct{,}
\]
we obtain:
\[
AVG
\le
(2\pi)^{-n/2}
\int_{p \in U}
\sqrt{
n! d\toric{n}
}
=
(2\pi)^{-n/2}
\int_{p \in U}
\sqrt{
n! d\toric{n}
}
\punct{.}
\]

Now we would like to use Cauchy-Schwartz again. This time,
the inner product is defined as:
\[
\langle \varphi , \psi \rangle \defeq \int_{p \in U} \varphi(p) \psi(p) dV
\punct{.}
\]

Hence,
\[
AVG
\le
(2\pi)^{-n/2}
\langle n! d\toric{n} , \mathbf{1} \rangle
\le
(2\pi)^{-n/2}
\| n! d\toric{n} \| \|\mathbf{1}\|
\punct{.}
\]

This time, $\|\mathbf{1}\|^2 = \lambda(U)$, so
we bound:

\begin{eqnarray*}
AVG
&\le&
(2\pi)^{-n/2}
\sqrt{\lambda(U)}
\sqrt{
\int_U n! d\toric{n}
}
\\
&\le&
(4\pi^2)^{-n/2}
\sqrt{\lambda(U)}
\sqrt{
\int_{(p,q)\in \toric{n}, p\in U} n! d\toric{n}
}. \text{\qed} 
\end{eqnarray*}

\subsection{The Proof of Theorem~\ref{unmixed:real}}
Let $\eps > 0$. As in the mixed case, we define:
\begin{eqnarray*}
\nu_{\mathbb R}(U,\eps) &\defeq&
\mathrm{Prob}_{f \in \mathcal F} 
\left[
\mymu(f;U) > \eps^{-1}
\right]
\\
&=&
\mathrm{Prob}_{f \in \mathcal F} 
\left[
\exists p \in U : \mathit{ev}(f;p)=0 \text{\ and \ }
d_{\mathbb P} (f, \Sigma_{p}) < \eps
\right]
\end{eqnarray*}
  where now $U \in \mathbb R^n$.

Let $V(\eps) \defeq \{ (f,p) \in \mathbb \mathcal F_{\mathbb R} \times U: 
\mathit{ev}(f;p)=0 \text{\ and \ }
d_{\mathbb P} (f, \Sigma_{p}) < \eps \}$. We also define $\pi: V(\eps)
\rightarrow \mathbb P(\mathcal F)$ to be the canonical projection 
mapping $F_\R\times U$ to $F_\R$ and set 
$\#_{V(\eps)}(f) \defeq \# \{ p \in U: (f,p) \in V(\eps) \}$.
Then,

\begin{eqnarray*}
\nu_{\mathbb R}(U,\eps)
&=&
\int_{f \in \mathcal F^{\mathbb R}}
\frac{e^{-\sum_i \|f^i\|^2/2}}{\sqrt{2\pi}^{\sum M_i}}
\chi_{\pi(V(\eps))} (f) \ 
d\mathcal F^{\mathbb R} 
\\
&\le&
\int_{f \in \mathcal F^{\mathbb R}}
\frac{e^{-\sum_i \|f^i\|^2/2}}{\sqrt{2\pi}^{\sum M_i}}
\#_{V(\eps)} 
d\mathcal F^{\mathbb R} 
\\
&\le&
\int_{p\in U \subset \mathbb R^n}
\int_{
\substack{
f \in \mathcal F^{\mathbb R}_{p} \\
d_{\mathbb P} (f, \Sigma_{p}) < \eps
}}
\frac{e^{-\sum_i \|f^i\|^2/2}}{\sqrt{2\pi}^{\sum M_i}}
\frac{1}{NJ(f;p)} 
d\mathcal F^{\mathbb R}_{p} 
\ dV_{\toric{n}}
\end{eqnarray*}

As before, we change coordinates in each fiber of
$\mathcal F^{\mathbb R}_{A}$ by
\[
f = f_{\romone} + f_{\romtwo} + f_{\romthree}
\]
with $f^i_{\romone}$ colinear to $v_{A}^T$,
$(f^i_{\romtwo})^T$ in the range of $Dv_{A}$,
and $f^i_{\romthree}$ othogonal to $f^{i}_{\romone}$
and $f^i_{\romtwo}$.
This coordinate system is dependent on $p+q\sqrt{-1}$.

In the new coordinate system, formula~\ref{det:cond:matrix}
splits as follows: 
\begin{eqnarray*}
\det \left(DG_{\scriptscriptstyle (p)} DG_{\scriptscriptstyle (p)}^H\right)^{-1/2}
dV_{\toric{n}} =
\hspace{-97pt} && \\
&=&
\left| \det 
\left[
\begin{matrix}
(f^1_{\romtwo})_1 & \hdots & (f^1_{\romtwo})_n\\ 
\vdots & & \vdots\\
(f^n_{\romtwo})_1 & \hdots & (f^n_{\romtwo})_n\\ 
\end{matrix}
\right]
\right|
\left| \det 
\left[
\begin{matrix}
({Dv_{A}}^{\romtwo})^1_1 & \hdots & ({Dv_{A}}^{\romtwo})^1_n  \\ 
\vdots & & \vdots \\
({Dv_{A}}^{\romtwo})^n_1 & \hdots & ({Dv_{A}}^{\romtwo})^n_n  \\ 
\end{matrix}
\right]
\right|
dV
\\
&=&
\left| \det 
\left[
\begin{matrix}
(f^1_{\romtwo})_1 & \hdots & (f^1_{\romtwo})_n\\ 
\vdots & & \vdots\\
(f^n_{\romtwo})_1 & \hdots & (f^n_{\romtwo})_n\\ 
\end{matrix}
\right]
\right|
\sqrt{\det Dv_A^H Dv_A}
\end{eqnarray*}

The integral $E(U)$ of $\sqrt{\det Dv_A Dv_A^H}$ is the expected number
of real roots on $U$, 
therefore 
\begin{multline*}
\nu_{\mathbb R}(U,\eps)
\le
E(U)
\int_{
\substack{
f_{\romtwo}+f_{\romthree} \in \mathcal F^{\mathbb R}_{p} \\
d_{\mathbb P} (f_{\romtwo}+f_{\romthree}, \Sigma_{p}) < \eps
}}
\frac{e^{-\sum_i \|f_{\romtwo}^i + f_{\romthree}^i\|^2/2}}{\sqrt{2\pi}^{\sum M_i}}
\cdot
\\
\cdot
\left| \det 
\left[
\begin{matrix}
(f^1_{\romtwo})_1 & \hdots & (f^1_{\romtwo})_n\\ 
\vdots & & \vdots\\
(f^n_{\romtwo})_1 & \hdots & (f^n_{\romtwo})_n\\ 
\end{matrix}
\right]
\right|
\ d\mathcal F^{\mathbb R}_{p} \ 
\punct{.}
\end{multline*}

In the new system of coordinates,
$\Sigma_{p}$ is defined by the 
equation:
\[
\det 
\left[
\begin{matrix}
(f^1_{\romtwo})_1 & \hdots & (f^1_{\romtwo})_n\\ 
\vdots & & \vdots\\
(f^n_{\romtwo})_1 & \hdots & (f^n_{\romtwo})_n\\ 
\end{matrix}
\right]
=
0
\punct{.}
\]

Since $\| f_{\romtwo} + f_{\romthree} \| \ge \| f_{\romtwo}\|$,
\[
d_{\mathbb P} (f_{\romtwo}+f_{\romthree}, \Sigma_{p}) < \eps
\Longrightarrow
d_{\mathbb P} (f_{\romtwo}, \Sigma_{p}) < \eps
\punct{.}
\]

This implies:
\begin{multline*}
\nu_{\mathbb R}(U,\eps)
\le
E(U) 
\int_{
\substack{
f_{\romtwo}+f_{\romthree} \in \mathbb \mathcal F^{\mathbb R}_{p} \\
d_{\mathbb P} (f_{\romtwo}, [\det = 0]) < \eps
}}
\frac{e^{-\sum_i \|f_{\romtwo}^i + f_{\romthree}^i\|^2/2}}{\sqrt{2\pi}^{\sum M_i}}
\cdot
\\
\cdot
\left| \det 
\left[
\begin{matrix}
(f^1_{\romtwo})_1 & \hdots & (f^1_{\romtwo})_n\\ 
\vdots & & \vdots\\
(f^n_{\romtwo})_1 & \hdots & (f^n_{\romtwo})_n\\ 
\end{matrix}
\right]
\right|
\ d\mathcal F^{\mathbb R}_{p} \ 
\punct{.}
\end{multline*}

We can integrate the $(\sum M_i - n - 1)$ variables $f_{\romthree}$ to
obtain:
\[
\nu_{\mathbb R}(U,\eps)
=
E(U)
\int_{\substack{
f_{\romtwo} \in \mathbb R^{n^2}\\
d_{\mathbb P} (f_{\romtwo}, [\det = 0]) < \eps
}}
\frac{e^{-\sum_i \|f_{\romtwo}^i\|^2/2}}{\sqrt{2\pi}^{n^2}}
\left| \det f_{\romtwo} \right|^2
\ d\mathbb R^{n^2} 
\punct{.}
\]

This is $E(U)$ times the probability
$\nu(n,\eps)$ for the linear case. 
\qed

\section{Acknowledgements}
Steve Smale provided valuable inspiration for us to extend 
the theory of \cite{BEZ1,BEZ2,BEZ3,BEZ4,BEZ5}
to sparse polynomial systems. He also provided
examples on how to eliminate the dependency upon unitary
invariance in the dense case.

  The paper by Gromov~\cite{GROMOV} was of foremost importance to
this research. To the best of our knowledge, ~\cite{GROMOV} is the only clear
exposition available of mixed volume in terms of a wedge of differential
forms. We thank Mike Shub for pointing out that reference and for
many suggestions.

  We would like to thank Jean-Pierre Dedieu for sharing his thoughts
with us on Newton iteration in Riemannian and quotient manifolds.  

  Also, we would like to thank Felipe Acker, 
Felipe Cucker, Alicia Dickenstein,
Ioannis Emiris, Askold Khovanskii, Eric Kostlan, T.Y. Li, Nelson
Maculan, Martin Sombra and Jorge P.\ Zubelli for their suggestions and support.

This paper was written while G.M. was visiting the Liu Bie Ju Center
for Mathematics at the City University of Hong Kong. He wishes to thank
City U for its generous support.

\appendix
\section{The Coarea Formula}
Here we give a short proof of the coarea formula,
in a version suitable to the setting of this paper. This means
we take all manifolds and functions smooth and avoid measure
theory as much as possible.

\begin{proposition} \label{coarea2}
\mbox{}\\
\begin{enumerate}
\item
  Let $X$ be a smooth Riemann manifold, of dimension $M$ and 
volume form $|dX|$. 
\item
  Let $Y$ be a smooth Riemann manifold, of
dimension $n$ and volume form $|dY|$. 
\item
  Let $U$ be an open set of $X$, and $F: U \rightarrow Y$
be a smooth map, such that $DF_x$ is
surjective for all $x$ in $U$. 
\item
  Let $\varphi: X \rightarrow \mathbb R^+$ be a 
smooth function with compact support contained in $U$. 
\end{enumerate}
Then for almost all $z \in F(U)$, $V_z \defeq F^{-1}(z)$ is a smooth 
Riemann manifold, and
\[
\int_X \varphi(x) NJ(F;x) |dX|
=
\int_{z \in Y} 
\int_{x \in V_z}
\varphi(x) |dV_z| |dY|
\]
where $|dV_z|$ is the volume element of $V_z$ and
$NJ(F,x) = \sqrt{\det DF_x^H DF_x}$ is the product 
of the singular values of $DF_x$. \qed 
\end{proposition}

 By the implicit function theorem, whenever $V_z$ is non-empty, 
it is a smooth $(N-n)$-dimensional Riemann submanifold
of $X$. By the same reason, 
$V := \{ (z,x): x \in V_z \}$
is also a smooth manifold. 

Let $\eta$ be the following $N$-form restricted to $V$:
\[
\eta = dY \wedge dV_z
\punct{.}
\]

This is {\bf not} the volume form of $V$. 
The proof of Proposition~\ref{coarea2} is divided into two
steps:

\begin{lemma}\label{coarea21}
\[
\int_V \varphi(x) |\eta| = \int_X \varphi(x) NJ(F;x) |dX|. 
\]
\end{lemma}

\begin{lemma}\label{coarea22}
\[
\int_V \varphi(x) |\eta| = \int_{z \in Y} \int_{x \in V_z}
\varphi(x) |dV_z| |dY|
\punct{.}
\]
\end{lemma}

\noindent {\bf Proof of Lemma~\ref{coarea21}:} 
We parametrize:
\[
\begin{array}{lrcl}
\psi: & X & \rightarrow & V \\
&x & \mapsto & (F(x),x)
\end{array}
\punct{.}
\]

Then,
\[
\int_V \varphi (x) |\eta| = \int_X (\varphi \circ \psi) (x) |\psi ^* \eta|
\punct{.}
\]

We can choose an orthonormal basis $u_1, \cdots, u_M$ of $T_xX$ such that
$u_{n+1}, \cdots, u_M \in \ker DF_x$. Then,
\[
D\psi (u_i) = 
\left\{
\begin{array}{ll}
 (DF_x u_i, u_i) & i=1, \cdots, n \\
 (0,u_i)         & i=n+1, \cdots, M
\end{array}
\right.
\punct{.}
\]
Thus,
\begin{eqnarray*}
| \psi^* \eta (u_1, \cdots, u_M) | &=&
| \eta (D \psi u_1, \cdots, D\psi u_M) | \\
&=&
|dY( DF_x u_1, \cdots, DF_x u_n )| 
\ |dV_z(u_{n+1}, \cdots, u_M)|
\\
&=& 
|\det DF_x |_{\ker DF_x^{\perp}}| \\
&=&
NJ(F,x) 
\end{eqnarray*}
and hence
\[
\int_V \varphi(x) |\eta| = \int_X \varphi(x) NJ(F;x) |dX|. \ \text{\qed} 
\]

\noindent {\bf Proof of Lemma~\ref{coarea22}:} 
We will prove this Lemma locally, and this 
implies the full Lemma through a standard
argument (partitions of unity in a compact
neighborhood of the support of $\varphi$). 
\medskip

Let $x_0,z_0$ be fixed. A small enough neighborhood of
$(x_0,z_0) \subset V_{z_0}$ admits a fibration
over $V_{z_0}$ 
by planes orthogonal to $\ker DF_{x_0}$.

We parametrize:
\[
\begin{array}{llcl}
\theta: & Y \times V_{z_0} & \rightarrow & V \\
& (z,x) & \mapsto & (z, \rho(x,z))
\end{array}
\punct{,}
\]

where $\rho(x,z)$ is the solution of $F(\rho)=z$ in the fiber
passing through $(z_0,x)$.
Remark that $\theta^* dY = dY$, and $\theta^* dV_z = \rho^* DV_z$.
Therefore,
\[
\theta^* (dY \wedge dV_z) = dY \wedge (\rho^* dV_{z})
\punct{.}
\]

Also, if one fixes $z$, then $\rho$ is a parametrization
$V_{z_0} \rightarrow V_z$. We have:
\begin{eqnarray*}
\int_V \varphi(x) |\eta| 
&=& 
\int_{Y \times V_{z_0}} \varphi(\rho(x,z)) |\theta^* \eta|
\\
&=& 
\int_{z \in Y} 
\left(
\int_{x \in V_{z_0}} \varphi(\rho(x,z) |\rho^* dV_{z}|
\right) |dY|
\\
&=&
\int_{z \in Y} 
\left(
\int_{x \in V_{z}} \varphi(x) |dV_{z}|
\right) |dY| \text{\qed} 
\end{eqnarray*}

The proposition below is essentially Theorem~3 p.~240 of~\cite{BCSS}.
However, we do not require our manifolds to be compact. We assume
all maps and manifolds are smooth, so that we can apply 
proposition~\ref{coarea2}.

\begin{proposition}\label{coarea1}\mbox{}\\
\vspace{-.5cm} 
\begin{enumerate}
\item \label{ca1} Let $X$ be a smooth $M$-dimensional manifold with volume element $|dX|$.
\item \label{ca2} Let $Y$ be a smooth $n$-dimensional manifold with volume element $|dY|$.
\item \label{ca3} Let $V$ be a smooth $M$-dimensional submanifold of $X \times Y$, and
     let $\pi_1: V \rightarrow X$ and $\pi_2: V \rightarrow Y$ be the 
canonical projections from $X\times Y$ to its factors. 
\item \label{ca4} Let $\Sigma'$ be the set of critical points of $\pi_1$,
     we assume that $\Sigma'$ has measure zero and that $\Sigma'$ is
     a manifold.
\item \label{ca5} We assume that $\pi_2$ is regular (all points in $\pi_2(V)$ are regular
     values).
\item \label{ca6} For any open set $U \subset V$, for any $x \in X$, 
     we write: $\#_U(x) \defeq \# \{ \pi_1^{-1}(x) \cap U \}$. We
     assume that $\int_{x \in X} \#_V (x) |dX|$ is finite. 
\end{enumerate}

     Then, for any open set $U \subset V$,
\[
\int_{x \in \pi_1(U)} \#_U (x) |dX| 
=
\int_{z \in Y}
\int_{\substack{x \in V_z\\ (x,z) \in U}} \frac{1}{\sqrt{ \det DG_x DG_x^H }} |dV_z| |dY|  
\]
  where $G$ is the implicit function for $(\hat x,G(\hat x)) \in V$ in a 
neighborhood of $(x,z) \in V \setminus \Sigma'$. \qed 
\end{proposition}

\noindent {\bf Proof:}
  Every $(x,z) \in U \setminus \Sigma'$ admits an open neighborhood such that
$\pi_1$ restricted to that neighborhood is a diffeomorphism. This defines
an open covering of $U \setminus \Sigma'$. Since $U \setminus \Sigma'$ is
locally compact, we can take a countable subcovering and define a partition of
unity $(\varphi_{\lambda})_{\lambda \in \Lambda}$ subordinated to that 
subcovering.

Also, if we fix a value of $z$, then $(\varphi_{\lambda})_{\lambda \in \Lambda}$
becomes a partition of unity for $\pi_1(\pi_1^{-1} (V_z) \cap U)$. Therefore,
\begin{eqnarray*}
\int_{x \in \pi_1(U)} \#_U (x) |dX| 
&=&
\sum_{\lambda \in \Lambda}
\int_{x,z \in \supp \varphi_{\lambda}} \varphi_{\lambda}(x,z) |dX| 
\\
&=&
\sum_{\lambda \in \Lambda}
\int_{z \in Y}
\int_{x,z \in \supp \varphi_{\lambda}} \frac{\varphi_{\lambda}(x,z)}{NJ(G,x)} |dX| 
\\
&=&
\int_{z \in Y}
\sum_{\lambda \in \Lambda}
\int_{x,z \in \supp \varphi_{\lambda}} \frac{\varphi_{\lambda}(x,z)}{NJ(G,x)} |dX| 
\\
&=&
\int_{z \in Y}
\int_{x \in V_z} \frac{1}{NJ(G,x)} |dX| 
\end{eqnarray*}

\noindent
where the second equality uses Proposition~\ref{coarea2} with $\varphi = \varphi_{\lambda} / NJ$.
Since $NJ = \sqrt{ \det DG_x DG_x^H }$, we are done.
\qed

\bibliographystyle{amsalpha}
\footnotesize

\providecommand{\bysame}{\leavevmode\hbox to3em{\hrulefill}\thinspace}
\providecommand{\MR}{\relax\ifhmode\unskip\space\fi MR }
\providecommand{\MRhref}[2]{%
  \href{http://www.ams.org/mathscinet-getitem?mr=#1}{#2}
}
\providecommand{\href}[2]{#2}

\end{document}